\documentclass[12pt]{article}

\usepackage{amscd,amsmath, amssymb, fancyhdr, url,epsfig, color, enumitem, tocloft}
\usepackage{dutchcal}
\usepackage{graphicx}

\usepackage[backref=page]{hyperref}
\renewcommand*{\backref}[1]{}
\renewcommand*{\backrefalt}[4]{%
	\ifcase #1 (Not cited.)%
	\or        (Cited on page~#2.)%
	\else      (Cited on pages~#2.)%
	\fi}

\hypersetup{
	colorlinks   = true,
	citecolor    = magenta,
	linkcolor    =blue,
	urlcolor     =magenta	
}

\newcommand{\version}{version 2.8,\ \ April 28, 2020}
\makeatletter
\def\x@arrow{\DOTSB\Relbar}
\def\xlongrightarrowfill@{\arrowfill@\relbar\relbar\longrightarrow}
\newcommand{\xlongrightarrow}[2][]{%
        \ext@arrow 0099\xlongrightarrowfill@{#1}{#2}}
\makeatother

\newcommand{\Lie}{\operatorname{Lie}}
\newcommand{\Jac}{\operatorname{Jac}}

\numberwithin{equation}{section}

\def\eqref#1{(\ref{#1})}
\newcommand{\goth}{\mathfrak}

\newcommand{\Z}{{\mathbb Z}}
\newcommand{\C}{{\mathbb C}}
\newcommand{\R}{{\mathbb R}}
\newcommand{\Q}{{\mathbb Q}}

\newcommand{\6}{\partial}
\def\1{\sqrt{-1}\:}
\newcommand{\restrict}[1]{{\left|_{{\phantom{|}\!\!}_{#1}}\right.}}
\newcommand{\cntrct}                
{\hspace{2pt}\raisebox{1pt}{\text{$\lrcorner$}}\hspace{2pt}}
\newcommand{\arrow}{{\:\longrightarrow\:}}

\newcommand{\calo}{{\cal O}}

\renewcommand{\bar}{\overline}
\renewcommand{\phi}{\varphi}
\renewcommand{\epsilon}{\varepsilon}
\renewcommand{\geq}{\geqslant}
\renewcommand{\leq}{\leqslant}

\newcommand{\im}{\operatorname{im}}

\newcommand{\Tot}{\operatorname{Tot}}

\newcommand{\ind}{\operatorname{\text{\sf ind}}}
\newcommand{\hor}{{\operatorname{\text{\sf hor}}}}
\renewcommand{\vert}{{\operatorname{\text{\sf vert}}}}

\newcommand{\Vol}{\operatorname{Vol}}
\newcommand{\Hom}{\operatorname{Hom}}
\newcommand{\Aut}{\operatorname{Aut}}
\newcommand{\Alb}{\operatorname{Alb}}

\newcommand{\coker}{\operatorname{coker}}

\newcommand{\rk}{\operatorname{rk}}

\newcounter{Mycounter}[section]
\newcounter{lemma}[section]
\renewcommand{\thelemma}{{Lemma \thesection.\arabic{lemma}}}
\newcommand{\lemma}{%
     \setcounter{lemma}{\value{Mycounter}}
     \refstepcounter{lemma}
     \stepcounter{Mycounter}
     {\noindent \bf \thelemma:\ }}

\newcounter{claim}[section]
\renewcommand{\theclaim}{{Claim \thesection.\arabic{claim}}}
\newcommand{\claim}{%
     \setcounter{claim}{\value{Mycounter}}
     \refstepcounter{claim}
     \stepcounter{Mycounter}
     {\noindent \bf \theclaim:\ }}

\newcounter{sublemma}[section]
\newcounter{corollary}[section]
\renewcommand{\thecorollary}{{Corollary \thesection.\arabic{corollary}}}
\newcommand{\corollary}{%
     \setcounter{corollary}{\value{Mycounter}}
     \refstepcounter{corollary}
     \stepcounter{Mycounter}
     {\noindent \bf \thecorollary:\ }}

\newcounter{theorem}[section]
\renewcommand{\thetheorem}{{Theorem \thesection.\arabic{theorem}}}
\newcommand{\theorem}{%
     \setcounter{theorem}{\value{Mycounter}}
     \refstepcounter{theorem}
     \stepcounter{Mycounter}
     {\noindent \bf \thetheorem:\ }}

\newcounter{conjecture}[section]
\newcounter{proposition}[section]
\renewcommand{\theproposition} {{Proposition \thesection.\arabic{proposition}}}
\newcommand{\proposition}{%
     \setcounter{proposition}{\value{Mycounter}}
     \refstepcounter{proposition}
     \stepcounter{Mycounter}
     {\noindent \bf \theproposition:\ }}

\newcounter{definition}[section]
\renewcommand{\thedefinition} {{Definition~\thesection.\arabic{definition}}}
\newcommand{\definition}{%
     \setcounter{definition}{\value{Mycounter}}
     \refstepcounter{definition}
     \stepcounter{Mycounter}
     {\noindent \bf \thedefinition:\ }}

\newcounter{example}[section]
\newcounter{remark}[section]
\renewcommand{\theremark}{{Remark \thesection.\arabic{remark}}}
\newcommand{\remark}{%
     \setcounter{remark}{\value{Mycounter}}
     \refstepcounter{remark}
     \stepcounter{Mycounter}
     {\noindent \bf \theremark:\ }}

\newcounter{problem}[section]
\newcounter{question}[section]
\makeatletter

\@addtoreset{equation}{section}
\@addtoreset{footnote}{section}
\makeatother

\def\blacksquare{\hbox{\vrule width 5pt height 5pt depth 0pt}}
\def\endproof{\blacksquare}

\newcommand{\proof}{{\bf Proof: \ }}
\newcommand{\pstep}{{\bf Proof. Step 1: \ }}

\begin{document}

\begin{center}
{\Large\bf  Classification of non-K\"ahler surfaces and\\[3mm]
locally conformally K\"ahler geometry}\\[5mm]
{\large
Liviu Ornea\footnote{Liviu Ornea and Victor Vuletescu are partially supported by a grant of Ministry of Research and Innovation, CNCS - UEFISCDI, project number PN-III-P4-ID-PCE-2016-0065, within PNCDI III.},  
Misha
Verbitsky\footnote{Misha Verbitsky is partially supported by the 
Russian Academic Excellence Project '5-100'' 
and CNPq - Process 313608/2017-2.\\[1mm]

}\\[4mm]
and Victor Vuletescu$^1$.
}

\end{center}

{\small
\hspace{0.15\linewidth}
\begin{minipage}[t]{0.7\linewidth}
{\bf Abstract} \\ 
Enriques-Kodaira classification treats
non-K\"ahler surfaces as a special case within the
Kodaira's framework. We prove the classification
results for non-K\"ahler complex surfaces not
relying on the machinery of Enriques-Kodaira classification.
We deduce the classification theorem of 
non-K\"ahler surfaces from the Buchdahl-Lamari theorem.
We also prove that all non-K\"ahler surfaces which are
not of class VII are locally conformally K\"ahler.
\end{minipage}
}


\tableofcontents


\section{Introduction}
\label{_Intro_Section_}

This paper appeared as a distillation of a lecture course
on complex surfaces given in 2008 and 2012 in Moscow Independent
University. The main reference on complex surfaces is the
great book by Barth, Hulek, Peters and Van de Ven 
\cite{_Barth_Peters_Van_de_Ven_}. This book offers
a powerful narrative, but (as it often happens with
great books) some plots are nested within more
plots, and it sometimes becomes hard to separate a 
particular strain from the polyphonic discourse.

For a forthcoming book on locally conformally K\"ahler 
(LCK) manifolds the first named two authors needed a classification of LCK structures
on (a posteriori, non-K\"ahler) surfaces. It was easier
(and more enlightening) to prove the non-K\"ahler part
of the Kodaira-Enriques classification directly along
with the classification of LCK structures.\footnote{For
a definition and an introduction to LCK structures, see 
Subsection \ref{_LCK_Subsection_}.}
 There are (almost) no new results of this paper, but most
of the proofs are different from those given in 
the literature, such as \cite{_Barth_Peters_Van_de_Ven_}. 

We tried to keep this expos\'e self-contained.
With the exception of two results from  Lamari's
paper \cite{_Lamari_}, we invoke only general notions
of complex algebraic geometry, found, for example,
in \cite{_Demailly:book_}.

The main focus of this paper is  elliptic non-K\"ahler
surfaces. We give a new proof that they are principal elliptic
bundles in the orbifold category (\ref{_elli_bundle_Theorem_}) - a result originally proven by Br\^inz\u anescu in \cite{_Brinzanescu:manuscripta_} ( see also \cite{_Brinzanescu:bundles_}) -, and locally conformally K\"ahler (\ref{_elli_then_Vaisman_Theorem_}) - a result due to Belgun, \cite{_Belgun_}. We prove that 
all these surfaces are Vaisman, giving a new proof
of Belgun's classification of Vaisman surfaces
(\cite{_Belgun_}).

We prove that all non-K\"ahler non-elliptic surfaces are 
of class VII (\ref{_elli_or_class_VII_Theorem_}). There is
enough good literature (\cite{_Dloussky_Oeljeklaus_Toma_},
\cite{_Dloussky:Kato_}, \cite{_Nakamura:towards_},
\cite{_Nakamura:surfaces_}, \cite{_Teleman:instantons_}) 
on class VII surfaces for us to 
give less attention to this case.

We state the known classification results up to the GSS conjecture
and give a new proof of Brunella's theorem on existence of LCK metrics
on Kato surfaces. Together with the GSS conjecture (still not fully proven)
this would imply that all non-K\"ahler complex surfaces are LCK,
with the exception of some of the Inoue surfaces
(\cite{_Belgun_}).

\subsection{Buchdahl-Lamari theorem}

In \cite{_Buchdahl:surfaces_,_Lamari_},
N. Buchdahl and A. Lamari have proven a result previously
known only from the Kodaira-Enriques classification of complex surfaces.

\hfill

\theorem\label{_B-L-intro_Theorem_}
Let $M$ be a compact complex surface. Then the first Betti number 
$b_1(M)$ is odd if and only if $M$ is non-K\"ahler.
\endproof

\hfill

Its direct proof, however, simplifies this 
classification significantly. In this paper we attempt
to recover most of the Kodaira-Enriques classification
for non-K\"ahler surfaces using the Buchdahl-Lamari theorem
and the following intermediate result 
(\ref{_Lamari-current-intro_Theorem_}), which was used
by Lamari to prove \ref{_B-L-intro_Theorem_}. Our proof is different
from the classical one, found in \cite{_Barth_Peters_Van_de_Ven_}, in a few aspects:
we do not rely on birational arguments and classification
of the elliptic surfaces due to Kodaira. Also, we 
aim to classify the locally conformally K\"ahler 
structures on complex surfaces.

\hfill

For an introduction to currents and their applications
in differential geometry, see \cite{_Demailly:book_}.
Recall that ``currents'' on $M$ are functionals on the space of 
differential forms on $M$ with compact support which are continuous in $C^\infty$-topology.
A differential form $\alpha$ defines a current $\tau \arrow \int_M \tau\wedge \alpha$.
This allows us to consider the differential forms as a subspace of  currents.

\hfill

\remark\label{_Current_cohomo_Remark_}
The usual operators and constructions of K\"ahler geometry 
(for example, $d, d^*, \6, \bar\6$, Laplacian, the Hodge decomposition)
extend from differential forms to currents in a natural way.
The corresponding cohomology (de Rham, Dolbeault, Bott-Chern) for
currents is equal to those of differential forms
(\cite{_Demailly:book_}).

\hfill

A (1,1)-form on a complex manifold $M$, $\dim_\C M=n$ is {\bf positive}
if is defined by a pseudo-Hermitian form with non-negative
eigenvalues, and {\bf strictly positive} when it is Hermitian. 
An $(n-1, n-1)$-form is {\bf positive} if it is
a product of $n-1$ positive forms. A (1,1)-current is called
{\bf positive} if it is non-negative on any positive
$(n-1, n-1)$-form.\footnote{For historical reasons, ``positivity''
for differential forms is understood in French sense:
0 is ``positive.'' We idly suggest the term ``French-positive'', 
to avoid confusion.} The notion of positivity for forms is 
compatible with that of currents (\cite{_Demailly:book_}).

\hfill

\theorem\label{_Lamari-current-intro_Theorem_}
Let $M$ be a compact complex non-K\"ahler surface. Then
there exists a non-zero positive, exact (1,1)-current $\Theta$ on $M$.

\proof \cite[Theorem 6.1]{_Lamari_}. \endproof

\hfill

\remark\label{_exact_pos_non-Kahler_Remark_}
The existence of an exact, non-zero positive (1,1)-current $\Theta$
immediately implies that the surface $M$ is non-K\"ahler.
Indeed, suppose that $M$ admits a K\"ahler form $\omega$.
Then $\int_M \omega\wedge \Theta > 0$, because $\omega\wedge \Theta$
is a non-zero measure on $M$, called {\bf the mass measure}
(\cite[Chapter III, Remark 1.15]{_Demailly:book_}). However, 
the strict inequality is impossible since $\Theta$ is exact.

\subsection{Locally conformally K\"ahler surfaces}
\label{_LCK_Subsection_}

Let us recall that a complex 
manifold $M$, $\dim_\C M >1$, is called {\bf locally conformally K\"ahler}
(LCK) if it admits a Hermitian form $\omega$ such that $d\omega = \theta \wedge \omega$,
where $\theta$ is a closed 1-form, see \cite{_Dragomir_Ornea_}. Then any cover $\tilde
M \stackrel{ \tau}{\arrow} M$ such that $\tau^*(\theta)=df$ is exact
is K\"ahler, with the K\"ahler form given by 
$\tilde \omega := e^{-f} \tau^*\omega$.  When $\theta$ is exact, $M$ is called
{\bf globally conformally K\"ahler}.  However, when $\theta$ is not
exact, and $M$ is compact, $M$ is non-K\"ahler
(\cite{_Vaisman:non-Kahler_}). If, in addition, there
exists a holomorphic conformal flow $\rho:\; \C \arrow
\Aut(M)$ acting on the K\"ahler cover $(\tilde M, \tilde
\omega)$ non-isometrically, $M$ is called {\bf Vaisman}. 

\hfill

\remark
Since $(\tilde M, \tilde \omega)$ is K\"ahler,
any LCK manifold can be obtained as a quotient
of a K\"ahler manifold by a discrete, proper action
of a group of holomorphic automorphisms, acting
on $\omega$ by holomorphic homotheties.\footnote{Notice that any conformal
holomorphic map of a connected K\"ahler manifold $\phi:\; (\tilde M,
\tilde \omega)\arrow (\tilde M, \tilde \omega)$, with
$\dim_\C \tilde M >1$, takes
$\tilde \omega$ to $f \tilde \omega$, with $df=0$
because $d\phi^*\tilde \omega = d(f\tilde \omega) =
df\wedge \omega=0$. Therefore, any conformal automorphism
of a K\"ahler manifold is a homothety.}
Thus, an LCK manifold can be defined
as a quotient of a K\"ahler manifold
$(\tilde M, \tilde \omega)$ by a
discrete, proper group acting by holomorphic
homotheties.

\hfill

\definition
Let $M$ be a compact complex surface with $b_1(M)=1$.
It is called {\bf a class VII surface} if its Kodaira
dimension is $\kappa(M)=-\infty$.

\hfill

From the Kodaira-Enriques classification
(\cite{_Barth_Peters_Van_de_Ven_}) it follows that
all minimal non-K\"ahler surfaces not of class VII
are elliptic. 
We prove this result (independently from the
rest of Kodaira-Enriques classification) in \ref{_elli_bundle_Theorem_}.
We prove that all elliptic surfaces are Vaisman when they are minimal.

\hfill

\theorem
Let $M$ be a compact, non-K\"ahler surface, 
which is not of class VII. Then $M$ admits an
LCK structure, and a Vaisman one if $M$ is
minimal.

\proof By \ref{_elli_or_class_VII_Theorem_}, $M$ is elliptic
when it is minimal, and by
\ref{_elli_then_Vaisman_Theorem_}, it is Vaisman. 
By \cite{_Tricerri_} and \cite{_Vuli:Bulletin_}, the blow-up at points 
preserves the LCK class, in particular a blow-up of an LCK
surface remains LCK. Therefore, a surface is LCK
if its minimal model is LCK, for example elliptic. 
\endproof

\hfill

For class VII surfaces, a complete classification
is not known, but it would follow from the so-called 
``Global Spherical Shell conjecture'' (GSS conjecture),
which claims that any minimal class VII surface $M$ with $b_2>0$
contains an open complex subvariety $U\subset M$
biholomorphic to a neighbourhood of the standard sphere 
$S^3 \subset \C^2$ such that $M \backslash U$ is connected.
Surfaces which satisfy this conjecture are called
{\bf Kato surfaces}. 

Brunella has shown that all Kato surfaces 
are LCK (\cite{_Brunella:Kato_} and Section \ref{_Brunella_proof_}).
Bogomolov's theorem on class VII
surfaces with $b_2(M)=0$ 
(see
\cite{_Bogomolov:VII_76_,_Bogomolov:VII_82_,
_Li_Yau_Zang:VII_,_Teleman:bogomolov_}) 
implies that they are
either Inoue surfaces or Hopf surfaces.
The Hopf surfaces are LCK (\cite{_GO:Hopf_surfaces_,_Belgun_,_OV:Shells_}), and
among the three classes of Inoue surfaces, two
are LCK, and the third contains a subclass which does not admit
an LCK structure (\cite{_Tricerri_,_Belgun_}). 

The modern proof of Bogomolov's classification theorem
(due to A. Teleman and Li-Yau-Zhang) is based on gauge
theory. Using gauge-theoretic methods, A. Teleman
was able to prove the GSS conjecture for minimal
class VII surfaces with $b_2=1$
(\cite{_Teleman:b2=1_}). Extending this
approach to $b_2>1$, Teleman was also able
to prove that any class VII manifold with
$b_2(M)=2$ contains a cycle of rational curves,
hence can be smoothly deformed to a blown-up Hopf surface
(\cite{_Teleman:instantons_}). 

Once the GSS conjecture is proven, this finishes the
classification of LCK surfaces. If it is true, all non-K\"ahler
surfaces are LCK, except a particular class  of Inoue surfaces,
which is non-LCK by results of F. Belgun (\cite{_Belgun_}).


\section{Cohomology of non-K\"ahler surfaces}
\label{_Cohomology_Section_}

\subsection{Bott-Chern cohomology of a surface }
\label{_BC_degree_Subsection_}

In this section, we repeat some of the 
arguments about the Bott-Chern cohomology previously given 
in \cite{_Teleman:cone_} and \cite{_Angella_Tomassini_V_}.
Recall that any compact complex manifold admits
a {\bf Gauduchon metric} in any conformal class of Hermitian
metrics (\cite{_Gauduchon_1984_}). By definition, a Gauduchon metric
on an $n$-dimensional complex manifold is a Hermitian metric
with Hermitian form $\omega$ satisfying $dd^c(\omega^{n-1})=0$,
where $d^c=I dI^{-1}$ is the twisted differential (here, 
as elsewhere, $I$ denotes the complex structure operator
extended to differential forms multiplicatively).

\hfill

\definition
The {\bf Bott-Chern cohomology group} of a complex manifold
is $H^{p,q}_{BC}(M):=\displaystyle\frac{\ker d\restrict{\Lambda^{p,q}(M)}}{\im dd^c}$.

\hfill

\remark
By the $dd^c$-lemma, 
on a compact K\"ahler manifold the Bott-Chern cohomology groups 
are equal to the de Rham and (hence) the Dolbeault
cohomology groups (\ref{_Current_cohomo_Remark_}). 
As usual, we denote the space of global
$(p,q)$-forms on $M$ by  $\Lambda^{p,q}(M)$.

It is well known (and not hard to see) that the
complex
\[
\Lambda^{p-1,q-1}(M)\stackrel{dd^c} \arrow \Lambda^{p,q}(M) \stackrel d\arrow \Lambda^{p+q+1}(M)
\]
is elliptic; see \cite[Proposition 5]{_Kodaira_Spencer:III_}. 
This implies that the Bott-Chern cohomology is  finite-dimensional
on any compact complex manifold.

\hfill

\theorem\label{_BC_Degree_Theorem_}
Let $M$ be a compact non-K\"ahler surface. Then the kernel of the
natural map $P:\;H^{1,1}_{BC}(M) \arrow H^2(M)$ is 1-dimensional.

\hfill

\pstep
Let $\omega$ be a Gauduchon metric on $M$.
Consider the differential operator $D:\; f\mapsto dd^c(f) \wedge \omega$ mapping
functions to 4-forms. Clearly, $D$ is elliptic and its index is the
same as  the index of the Laplacian: $\ind D = \ind \Delta=0$, hence $\dim \ker D=\dim \coker D$. 
The Hopf maximum principle implies that $\ker D$ only
contains constants, hence $\coker D$ is 1-dimensional.
However,
$\int_M D(f) = \int_M dd^c(f) \wedge \omega = \int_M fdd^c \omega=0$.
This implies that a 4-form $\kappa$ belongs to $\im D$ if and only if $\int_M\kappa=0$.

\hfill

{\bf Step 2:}
Let $\alpha$ be a closed $(1,1)$-form. Define {\bf the degree}
$\deg_\omega \alpha:= \int_M \omega \wedge\alpha$.
Since $\int_M dd^cf \wedge \omega=0$, this defines a map $\deg_\omega:\; H^{1,1}_{BC}(M,\R)\arrow \R$.
Given a closed $(1,1)$-form $\alpha$ of degree 0, 
the form $\alpha':=\alpha- dd^c(D^{-1}(\alpha\wedge \omega))$ 
satisfies $\alpha'\wedge \omega=0$, in other words, it is an $\omega$-primitive
(1,1)-form. For $\omega$-primitive forms, one has $\alpha'\wedge \alpha' = -|\alpha'|^2\omega\wedge\omega$,
giving
\begin{equation} \label{_primitive_square_Equation_}
\int_M \alpha'\wedge \alpha'= -\|\alpha'\|^2_\omega
\end{equation}
which is impossible when $\alpha'$ is a non-zero vector in $\ker P$,
because $\alpha'$ is exact. 
Therefore, any vector of zero degree in $\ker P\subset H^{1,1}_{BC}(M,\R)$ vanishes. This implies that
any two vectors in $\ker P$ are proportional.
\endproof

\subsection{First cohomology of non-K\"ahler surfaces}

The Dolbeault cohomology of a complex manifold $M$ 
is denoted by $H^{p,q}(M)$. Clearly, $H^{p,0}(M)$
coincides with the space of holomorphic $p$-forms on $M$.
On a compact K\"ahler manifold, holomorphic
$p$-forms are closed because they are harmonic.
On non-K\"ahler manifolds, this is generally false.
However, this is true on compact complex surfaces.

\hfill

\lemma\label{_holo_closed_Lemma_}
All holomorphic 1-forms on a compact complex surface are closed.

\hfill

\proof
Let $\alpha\in \Lambda^{1,0}(M)$ be a holomorphic 1-form. 
Then $\bar\6\alpha=0$, because it is holomorphic, and
by the same reason $d\alpha$ is a holomorphic, exact (2,0)-form. 
Then $d\alpha\wedge d\bar\alpha$ is a positive (2,2)-form,
giving $0=\int_M d\alpha\wedge d\bar\alpha =\|d\alpha\|^2$.
Then $d\alpha=0$, and $\alpha$ is closed.
\endproof

\hfill

\theorem\label{_H^1_odd_Theorem_}
Let $M$ be a non-K\"ahler manifold and
$\Theta$ a non-zero exact positive (1,1)-current, which exists by \ref{_Lamari-current-intro_Theorem_}.
Let $d\alpha = \beta$ be a real (1,1)-form
in the same Bott-Chern cohomology class. 
Then:
\begin{description}
\item[(i)] Denote the space of holomorphic 1-forms on $M$ by
$H^{1,0}(M)$, and let $\overline{H^{1,0}(M)}$ 
denote its complex conjugate (the space of antiholomorphic forms). 
By \ref{_holo_closed_Lemma_}, holomorphic 1-forms are closed, therefore,
there is a natural map $H^{1,0}(M)\oplus \overline{H^{1,0}(M)} \arrow H^1(M, \C)$ 
mapping the sum of a holomorphic and an antiholomorphic form to its
cohomology class. We claim that this map is injective, and,
moreover, for appropriate choice of $\alpha$ 
the form $\theta:= \alpha^{1,0} - \alpha^{0,1}$ is closed
and satisfies 
\begin{equation}\label{_Hodge_on_H^1_surface_Equation_}
H^1(M) = H^{1,0}(M)\oplus \overline{H^{1,0}(M)} \oplus
\langle[\theta]\rangle=H^1(M),
\end{equation}
that is, $H^1(M)$ is generated by cohomology classes
of holomorphic and antiholomorphic forms and $\theta$.
\item[(ii)] Since all holomorphic forms are closed,
the antiholomorphic forms are $\bar\6$-closed and have Dolbeault classes.
This gives a natural embedding $\overline{H^{1,0}(M)}\hookrightarrow H^{0,1}(M)$.
We claim that $H^{0,1}(M)$ is generated by $\overline{H^{1,0}(M)}$ and
the Dolbeault class $[\theta^{0,1}]$, giving
$H^{0,1}(M)= \overline{H^{1,0}(M)}\oplus \langle [\theta^{0,1}] \rangle$.
\end{description}

{\bf Proof of (i):} 
A non-zero linear combination of holomorphic and antiholomorphic 
forms is closed and never exact. Indeed, if $df$ is
a linear  combination of holomorphic and antiholomorphic 
forms, then $d^c df=0$, hence $f$ is a globally 
defined harmonic function on a compact manifold.
Such a function has to be constant by maximum principle.
Therefore,  the natural map 
$H^{1,0}(M)\oplus \overline{H^{1,0}(M)} \stackrel \kappa \arrow H^1(M,\C)$ 
is injective. To prove \ref{_H^1_odd_Theorem_} (i) it remains to show that
its image has codimension 1, and prove that $H^1(M,\C)$ 
is generated by the class of  $\theta$ and the image of $\kappa$. 

Clearly, the kernel of the natural map $\6:\; H^{0,1}(M)\arrow H^{1,1}_{BC}(M)$
coincides with $\overline{H^{1,0}(M)}$. Its image is contained in 
the one-dimensional kernel of the map $P:\;H^{1,1}_{BC}(M) \arrow H^2(M)$ 
(\ref{_BC_Degree_Theorem_}). Therefore, the image of
$\kappa$ has codimension 
at most 1. The class $\6 \alpha^{0,1}\in H^{1,1}_{BC}(M)$ 
has non-zero degree, because
$\6 \alpha^{0,1}+ \overline{\6 \alpha^{0,1}}=\beta$,
hence $\6 \alpha^{0,1}$ generates $\ker P= \langle [\beta] \rangle =
\langle [\Theta]\rangle$.

Replacing $\alpha^{0,1}$ by a Dolbeault cohomologous class
if necessary, we can always assume that $\beta=\6 \alpha^{0,1}$
and $\bar\6\alpha^{0,1}=0$. Since
$d\alpha$ is a real (1,1)-form, $\theta:= \alpha^{1,0} - \alpha^{0,1}$ is closed.
Since its (0,1)-part has the same Dolbeault
cohomology class as $-\alpha^{0,1}$,
the form $\theta$ is not exact. Therefore, it generates the cokernel of 
the embedding $H^{1,0}(M)\oplus \overline{H^{1,0}(M)} \stackrel \kappa
\arrow H^1(M,\C)$. 
We proved \eqref{_Hodge_on_H^1_surface_Equation_}.

\hfill

{\bf Proof of (ii):}
Consider the map $\6:\; H^{0,1}(M)\arrow H^{1,1}_{BC}(M)$.
By Step 1, $\6([\theta^{0,1}])$ is non-zero in
$H^{1,1}_{BC}(M)$.
Since $\6:\; H^{0,1}(M)\arrow H^{1,1}_{BC}(M)$
vanishes on antiholomorphic forms, 
the Dolbeault cohomology  class of $\theta^{0,1}$
does not belong to $\overline{H^{1,0}(M)}$,
giving an injective map 
\begin{equation}\label{_holo_to_Dolbeault_Equation_}
\overline{H^{1,0}(M)}\oplus \langle [\theta^{0,1}]
\rangle\arrow H^{0,1}(M).
\end{equation}
Now, the kernel of $\6:\; H^{0,1}(M)\arrow H^{1,1}_{BC}(M)$
is generated by antiholomorphic forms.
This gives an exact sequence
\begin{equation}\label{_exact_Dolbeault_Equation_}
0\arrow \overline{H^{1,0}(M)}\arrow H^{0,1}(M)\arrow
H^{1,1}_{BC}(M) 
\end{equation}
The image of $\6:\; H^{0,1}(M)\arrow H^{1,1}_{BC}(M)$
is at most 1-dimensional by \ref{_BC_Degree_Theorem_},
hence it is generated by $\6([\theta^{0,1}])$.
Using \eqref{_exact_Dolbeault_Equation_}, 
we obtain that the injective map 
\eqref{_holo_to_Dolbeault_Equation_}
is actually surjective. We proved 
\ref{_H^1_odd_Theorem_} (ii).
\endproof

\subsection{Second cohomology of non-K\"ahler surfaces}

\claim\label{_H^2,0_Claim_} 
Let  $M$ be a compact complex surface. Then: 
\[ 
   \overline{H^{0,2}(M)}= H^{2,0}(M) = H^0(\Omega^2(M)).
\]
\proof
By Serre's duality, $H^{0,2}(M) = H^0(K_M)^*$, which has the same
dimension as $H^0(\Omega^2(M)) = H^0(K_M)=H^{2,0}(M)$.
The natural map $R:\;H^{2,0}(M) \arrow \overline{H^{0,2}(M)}$ 
is injective, because its kernel is formed by $\6$-exact holomorphic
forms $\alpha= \6\beta$, but for such $\alpha$ one has 
$0=\int_Md\beta\wedge \bar \alpha = \int_M\alpha\wedge \bar\alpha= \|\alpha\|^2$,
which is impossible unless $\alpha=0$. As shown above, these spaces
have the same dimension, hence $R$ is an isomorphism.
\endproof

\hfill

\corollary 
(\cite[Theorem IV.2.8]{_Barth_Peters_Van_de_Ven_})
The Hodge-de Rham-Fr\"olicher spectral sequence of a compact complex
surface degenerates in $E_1$. 

\hfill

\proof This argument is standard: the Hodge-de Rham-Fr\"olicher 
degenerates if the dimension does not drop.

We have proven the degeneration 
for $H^1(M)$ in \ref{_H^1_odd_Theorem_}.
Serre's and Poincar\'e duality give the degeneration of this spectral sequence
in the page $E_1$ for the $H^3(M)$ term. 
Now it has to degenerate in 
the page $E_1$ in $H^2(M)$, because for each instance
where $d_k\neq 0$, with $d_k:\; E^{p,q}_{k-1}\arrow E^{p-k+1,q+k}_{k-1}$
non-zero, the group $E^{p,q}_{k-1}$ is replaced by $\ker d_k$ and
$E^{p-k+1,q+k}_{k-1}$ by $\coker d_k$. Then the total dimension
of the space $\bigoplus_{p+q=d}E^{p,q}_{k-1}$ and $\bigoplus_{p+q=d}E^{p-k+1,q+k}_{k-1}$
would decrease. However, for $d=1$ and $d=3$ it is already
minimal, hence $d_k=0$ for all $k\geq 2$.
\endproof

\hfill

\remark
We have shown that $H^2(M, \C)  = H^{2,0}(M) \oplus H^{1,1}(M) \oplus H^{0,2}(M)$
and $H^{2,0}(M)= \overline {H^{0,2}(M)}$, same as for K\"ahler surfaces
(\ref{_H^2,0_Claim_}). 

\hfill

\lemma\label{_BC_to_H^11_Lemma_}
Let $M$ be a compact complex surface. Then the natural
map $H^{1,1}_{BC}(M) \arrow H^{1,1}(M)$ 
(from Bott-Chern cohomology to Dolbeault cohomology)
is surjective.

\hfill

\proof
Let $\alpha\in \Lambda^{1,1}(M)$ be a $\bar\6$-closed form representing
a Dolbeault cohomology class $[\alpha]\in H^{1,1}(M)$.
Clearly, $\alpha$ can be represented by a closed (1,1)-form
if and only if the Bott-Chern class $[\beta]\in
H^{2,1}_{BC}(M)$ of $\beta:=\6 \alpha$ 
vanishes. Indeed, if $\beta=\6\bar\6\eta$, then
$\alpha - \6\eta$ is $\bar\6$-closed.

Since the Hodge - de Rham - Fr\"olicher 
spectral sequence degenerates in
$E_1$, and $\beta$ is 
$\6$-exact and $\bar\6$-closed,
one has $\beta= \bar\6\gamma$, for some $\gamma\in \Lambda^{2,0}(M)$.
However, the cokernel of the map
$\6:\; \Lambda^{1,0}(M)\arrow \Lambda^{2,0}(M)$
is generated by holomorphic 2-forms, and they are all
closed (\ref{_H^2,0_Claim_} ),
hence $\bar\6(\Lambda^{2,0}(M)) = \bar\6\6(\Lambda^{1,0}(M))$,
and the Bott-Chern class of $\beta\in \bar\6(\Lambda^{2,0}(M))$ vanishes.
\endproof

\hfill

\proposition\label{_neg_def_h^11_Proposition_}
(\cite[Theorem IV.2.14]{_Barth_Peters_Van_de_Ven_})
Let $M$ be a compact non-K\"ahler surface. Then the
intersection form on the image of $H^{1,1}_{BC}(M,\R)$
in de Rham cohomology is negative definite.

\hfill

\proof 
Fix a Gauduchon metric $\omega$ on $M$.
Let $[\alpha]\in H^{1,1}(M)$ be a cohomology class.
By \ref{_BC_to_H^11_Lemma_}, we can represent $[\alpha]$
by a closed (1,1)-form $\alpha$. Consider the degree
functional $\deg_\omega:\; H^{1,1}_{BC}(M,\R) \arrow \R$
defined in Subsection \ref{_BC_degree_Subsection_}.
Since $\deg_\omega(\Theta)>0$ for an exact (1,1)-current $\Theta$,
any cohomology class $[\alpha]\in H^{1,1}(M)\subset H^2(M)$
can be represented by a closed (1,1)-form $\alpha$ with
$\deg_\omega\alpha=0$. Acting as in the proof of
\ref{_BC_Degree_Theorem_}, we find 
$f\in C^\infty(M)$ such that $\alpha-dd^c f$ is primitive.
Replacing $\alpha$ by $\alpha-dd^c f$, we
obtain $\int_M \alpha\wedge \alpha= -\|\alpha\|^2_\omega<0$
as in \eqref{_primitive_square_Equation_}.
\endproof

\hfill

\remark
In \cite[Theorem 2.37]{_Brinzanescu:bundles_} 
(see also \cite{_Brinzanescu_Flondor:1_, _Brinzanescu_Flondor:3_}) it was shown that
the space $H^{1,1}(M) \cap H^2(M, \Q)$ (the rational Neron-Severi group) 
of a complex surface with algebraic
dimension 0 has negative definite intersection form.
This important result also follows from \ref{_neg_def_h^11_Proposition_}




\subsection{Vanishing of products of holomorphic 1-forms}

\ref{_neg_def_h^11_Proposition_} has an interesting corollary.

\hfill

\proposition \label{_prod_hol_1_forms_Proposition_}
Let $M$ be a  non-K\"ahler, compact, complex surface. Then for any 
holomorphic  1-forms $\alpha, \beta$, the product
$\alpha\wedge \beta$ vanishes, and the product $\alpha\wedge \bar\beta$ is exact.

\hfill

\pstep
Let $\alpha, \beta$ be holomorphic 1-forms.
Then $\alpha\wedge \beta$ is a holomorphic 2-form, with
$\int_M \alpha\wedge \beta\wedge \bar \alpha\wedge \bar \beta>0$
unless $\alpha\wedge \beta=0$. 

\hfill

{\bf Step 2:}
The intersection form on $H^{1,1}(M)$ is negative definite by
\ref{_neg_def_h^11_Proposition_}. 
Therefore, $\int_M \eta\wedge \bar\eta=0$ 
implies $\eta=0$ for any $\eta\in H^{1,1}(M)$.  
Take $\eta=\alpha\wedge \bar \alpha$
and $\rho=\beta\wedge \bar \beta$. Then clearly 
$\int_M \eta\wedge \bar\eta=\int_M \rho\wedge \bar\rho=0$,
which implies 
$\int_M \alpha\wedge \beta\wedge \bar \alpha\wedge \bar \beta=0$,
hence $\alpha\wedge \beta=0$ (Step 1).

\hfill

{\bf Step 3:}
To see that $\eta:=\alpha\wedge \bar\beta$ is exact,
we use \ref{_neg_def_h^11_Proposition_}
again. Unless $\mu:=\alpha\wedge \bar \beta$ is exact, one would have
$\int_M \mu\wedge \bar\mu<0$,
but $\int_M \mu\wedge \bar\mu=0$ 
as we have already shown in Step 2. Therefore, $\mu$ is exact.
\endproof

\subsection{Structure of multiplication in 
$H^1(M)$ for non-K\"ahler surface without curves}

Further on, we shall need the following lemma.

\hfill

\lemma\label{_A_defi_Lemma_}
Let $M$ be a compact non-K\"ahler surface,
and $\Theta$ an exact real (1,1)-form. Then there exists
a closed form $\theta \in \Lambda^1(M)$  such that 
$ d^c(\theta)=\Theta$. 

\hfill

\proof
Without restricting generality, we can assume that $\Theta$ 
is an exact real (1,1)-form representing a non-zero class in
$H^{1,1}_{BC}(M)$. By \ref{_H^1_odd_Theorem_} (ii), 
the cohomology class of $\Theta$ belongs to the
image of the natural map $\6:\; H^{0,1}(M) \arrow H^{1,1}_{BC}(M)$,
and can be expressed as $\6(\tilde \theta^{0,1})$, where
$\tilde\theta$ is a closed form. Therefore,
$\Theta- d^c(\tilde \theta)= d^cdf$  for some function $f\in
C^\infty M$. Take $\theta:= \tilde \theta+ d f$.
Then $\theta$ is a closed 1-form such that $d^c\theta=\Theta$.
\endproof

\hfill

In this subsection we prove the following structure theorem
for the multiplication in $H^1(M)$. It is a posteriori true
in the general situation, but we need it only for surfaces without
complex curves.

\hfill

\theorem \label{_multiplica_H^1_Theorem_}
Let $M$ be a compact non-K\"ahler surface
without curves, and $\theta \in \Lambda^1(M)$
a closed 1-form such that $d^c(\theta)$ is non-zero in $H^{1,1}_{BC}(M, \R)$
(\ref{_A_defi_Lemma_}).  Denote by $W \subset \Lambda^1(M)$
the subspace generated by holomorphic and antiholomorphic forms.
We identify $W$ with its image in 
$H^1(M)= W \oplus \langle \theta \rangle$ (\ref{_H^1_odd_Theorem_}).
Then 
\begin{description}
\item[(i)] The multiplication $W \wedge W \arrow H^2(M)$ vanishes.
The multiplication $W \wedge \theta \arrow H^2(M)$ is injective.
\item[(ii)] Let $\theta^c:= I(\theta)$.
Then for any non-zero $x \in W$, the form $x\wedge \theta^c$ is 
closed and represents a non-zero element of $H^2(M)$.
\item[(iii)] The Poincar\'e pairing on the spaces 
$W \wedge \theta$ and $W \wedge \theta^c\subset H^2(M)$ vanishes. The 
Poincar\'e pairing between these two subspaces is non-degenerate.
\end{description}

\pstep Multiplication $W \wedge W \arrow H^2(M)$ vanishes
by \ref{_prod_hol_1_forms_Proposition_}.
Let us prove that for any non-zero $x\in W$,
the form $x\wedge \theta^c$ is 
closed and represents a non-zero element of $H^2(M)$.

Without restricting the generality, we may assume that
$W\neq 0$. For any $x\in W$, the closed (1,1)-form $x\wedge\bar x$ 
is homologous to zero
(\ref{_prod_hol_1_forms_Proposition_}). 
Indeed, one has $(x\wedge\bar x)^2=0$,
and the Poincare pairing on the image of $H^{1,1}_{BC}(M)$ 
in the de Rham cohomology is negative definite 
by \ref{_neg_def_h^11_Proposition_}.

For any non-zero holomorphic form $x$, the form
$\Theta:=x\wedge \bar x$ is positive, non-zero and exact.
We fix such $\Theta$, and fix a closed form 
$\theta$ such that $d^c\theta=\Theta$ as in
\ref{_H^1_odd_Theorem_}  (ii). Choose any $y\in W$.
Then $d(y\wedge \theta^c) = y \wedge x\wedge \bar x=0$
(\ref{_prod_hol_1_forms_Proposition_}) hence
$y\wedge \theta^c$ is closed. 

\hfill

{\bf Step 2:}
The injectivity of $W \xlongrightarrow{x\wedge \theta} H^2(M)$
and $W \xlongrightarrow{x\wedge \theta^c} H^2(M)$ would follow if
we prove that the formula 
\[
x, y \arrow \int_M x\wedge y \wedge \theta\wedge \theta^c,
\ \ x, y\in W
\]
defines a non-degenerate pairing on $W$,
these two spaces is non-degenerate, hence (i) and (ii) follows from (iii).

The Poincar\'e pairing on the images $W\wedge \theta$ of 
$W\wedge\theta^c$ in $H^2(M)$ 
vanishes because $\theta\wedge\theta = \theta^c \wedge\theta^c=0$.
To prove that the Poincar\'e pairing between $W\wedge \theta$
and $W\wedge \theta^c$ is non-degenerate, take a holomorphic 
1-form $x\in W$. To finish the proof of \ref{_multiplica_H^1_Theorem_} (iii) it 
would suffice to show that the integral 
$\int_M \1 x\wedge \bar x \wedge \theta \wedge \theta^c$
is positive.  The form $\1 x\wedge \bar x \wedge \theta
\wedge \theta^c$ is positive, because it is
a product of a $(2,0)$-form and its complex conjugate.
It is  non-zero if $x$ is not
proportional to the (1,0)-part of $\theta$, denoted 
$\theta^{1,0}$. 

\hfill

{\bf Step 3:}
It remains to show that $\theta^{1,0}$
is not proportional to a holomorphic form; this is
where we use that $M$ has no holomorphic curves.
In this case the zero set of any 1-form $x\in W$ is 
a finite set.

Suppose that $x\in W$ is holomorphic
and proportional to $\theta^{1,0}$.
Then there exists a smooth function $\alpha$
defined outside of the zero set $S$ of $x$
such that $\theta^{1,0}= \alpha x$.
Without restricting generality, we may assume that
$d\theta^{1,0}=x \wedge \bar x$ (Step 1). Then 
$\bar \6\alpha = \bar x$, which gives $dd^c \alpha=0$
outside of $S$. Then $\alpha$ is locally
the real part of a holomorphic function;
using the Hartogs extension theorem, we obtain that
$\alpha$ is smooth and defined globally on $M$. Then $\alpha=0$
by maximum principle, because $dd^c$ is elliptic. 
\endproof

\hfill

\corollary\label{_W_rational_Corollary_}
Let $M$ be a compact surface
without curves, and $W \subset H^1(M, \C)$
the subspace generated by holomorphic and antiholomorphic forms.
Then $W$ is a rational subspace, that is, there
exists a subspace $W_\Q \subset H^1(M, \Q)$
such that $W = W_\Q \otimes_\Q \C$.

\hfill

\proof
For $x\in H^1(M)$, denote by $L_x:\; H^1(M)\arrow H^2(M)$
the map $y \arrow x\wedge y$. Without restricting
generality, we may assume $W\neq 0$; then $\dim W\geq
2.$ By \ref{_multiplica_H^1_Theorem_}, $\rk L_\theta= \dim
W\geq 2$
and $\rk L_x =1$ for $x\in W$. Therefore  
$W$ is the space of all $x\in  H^1(M)$ such that the
$\rk L_x=1$. Since the multiplication in
$H^*(M)$ is defined over $\Q$, the subspace $W\subset H^1(M)$
is rational.
\endproof

\hfill

\remark
\ref{_W_rational_Corollary_}
is true for all surfaces. We leave its proof
as an exercise to the reader.


\section{Barlet spaces}
\label{_Barlet_Section_}

Barlet spaces are spaces of cycles, that is, 
closed complex analytic subvarieties
of given dimension in a given complex manifold
with multiplicities (positive integers)
assigned to their irreducible components.
They are similar but distinct from the Douady spaces, which are
spaces of closed complex analytic subspaces (possibly with nilpotents
in the structure sheaf). 

For more details on Barlet spaces and their properties, please see the book
\cite{_Barlet_Magnusson:book_}.

Let $M$ be a metric space.
The (Gromov-)Hausdorff metric on 
the set ${\cal C}$ of closed subsets of $M$
is defined as follows: $d(X,Y)$ is infimum
of all $\epsilon$ such that $X$ belongs to
$\epsilon$-neighbourhood of $Y$, and
$Y$ belongs to $\epsilon$-neighbourhood of $X$.
When $M$ is compact, the corresponding topology
on ${\cal C}$ is independent from the choice
of metric on $M$ as long as topology of $M$ 
remains the same. It is called 
{\bf the (Gromov-)Hausdorff topology}.

It is important for the sequel that the
Barlet space of cycles is a reduced complex analytic space, with
the topology which is compatible with the (Gromov-)Hausdorff
topology on the set of all closed subvarieties.

For each $k \geq 0$, the Barlet space ${\goth B}_k(M)$ of $k$-cycles
on a manifold $M$ comes equipped with a universal family of cycles,
and their supports form a closed complex-analytic subvariety  ${\goth B}^m_k(M) \subset
M \times {\goth B}(M)$, the ``marked Barlet space''
of pairs ``complex analytic cycle and
a point in its support''. The
forgetful map ${\goth B}^m_k(M)\arrow {\goth B}_k(M)$ is equidimensional of
relative dimension $k$, and the forgetful map $\Psi:{\goth B}^m_k(K)
\to M$ is complex-analytic. In particular, for any compact irreducible
component $Z$ in a Barlet space, and its 
marked counterpart $Z^m$, the image $\Psi(Z^m)$ is a complex
analytic subvariety by Remmert's proper mapping
theorem. Geometrically,
this means that for any compact complex-analytic family $Z$ of
cycles in $M$, the union of supports of all these cycles
is a complex subvariety in $M$.

Further on, we shall use the following result, which
was proven for curves in \cite[Corollary 2.19]{_Verbitsky:twistor_},
and for divisors in \cite{_Barlet:divisors_}.

\hfill

\theorem\label{_compact_Barlet_Theorem_}
Let $(M,I)$ be a compact complex manifold, and
$\omega$ a Hermitian form. Assume that
$dd^c(\omega^k)=0$ for some integer $k>0$. Then any connected component
$Z$ of the Barlet space ${\mathfrak B}_k(M)$ 
of $k$-cycles in $M$ is compact.

\hfill

\proof
By Bishop's theorem (\cite{_Bishop:conditions_}),
a (Gromov-)Hausdorff limit of a family of compact complex
subvarieties is complex if its Hermitian volume
stays bounded. Since the space of closed
subvarieties in $M$ with (Gromov-)Hausdorff topology is compact,
this implies that the set of closed 
compact $k$-dimensional subvarieties
with volume bounded by a constant $C\in \R$
is compact.

Therefore, to prove compactness it would suffice to show
that the volume $\Vol(S)$ is constant as a function of 
$[S]\in Z$.

Let $X$ denote the marked family associated with $Z$,
and $\pi_M:\; X \arrow M$, $\pi_X:\; X \arrow Z$
the forgetful maps. Then the volume function
$\Vol:\; Z\arrow \R^{>0}$ can be expressed as
$\Vol = (\pi_X)_*\pi_M^* \omega^k$, where
$(\pi_X)_*$ denotes the pushforward of a
differential form (generally speaking, the
pushforward is not a form, but it is well
defined as a current).

Let $k$ be dimension of the cycles parametrized by $Z$.
Since pullback and pushforward of differential forms
commute with $d$, $d^c$, this gives $dd^c
\Vol=(\pi_X)_*\pi_M^* (dd^c\omega^k)$ (see e.g.
\cite[(8.12)]{_NHYM_}, \cite[Theorem 2.10]{_Verbitsky:S^6_}
or \cite[Proposition 1.9]{_Ivashkovich:Annals_}).
Therefore, $\Vol$ satisfies $dd^c(\Vol)=0$
whenever $dd^c(\omega^k)=0$.

Functions which satisfy $dd^c f=0$ are called
{\bf pluriharmonic}. Using local $dd^c$-lemma, it
is easy to see that any pluriharmonic function
is locally a sum of a holomorphic and antiholomorphic
function.

By Bishop's theorem, the set
$\Vol^{-1}(]-\infty, C])$ is compact for all $C\in \R$,
hence  $-\Vol$ has a maximum somewhere in $X$.
However, a pluriharmonic function which admits a 
maximum  is necessarily
constant by E. Hopf's strong maximum principle.
Therefore, $\Vol$ is constant on each
connected component of the Barlet space.
Now, each of these components is compact by Bishop.
\endproof

\hfill

Applying this result to the
space of curves on a complex surface
and using the Gauduchon form $\omega$
(Subsection \ref{_BC_degree_Subsection_}),
we obtain the following useful corollary.

\hfill

\corollary\label{_Barlet_on_surf_Corollary_}
Let $M$ be a compact complex surface,
${\goth B}_1(M)$ the Barlet space of 1-cycles on $M$
and $Z$ its connected component. Then $Z$ is compact.
\endproof


\section{Non-K\"ahler elliptic surfaces}
\label{_elliptic_Section_}

\subsection{The Gauss-Manin connection}

For the sequel, we recall some basic facts about the
Gauss-Manin connection. We may refer to \cite{_Griffiths:transcendental_} 
or \cite{_Voisin-Hodge_}. Let $\pi:\; M \arrow B$ be
a smooth, proper map of smooth manifolds. By Ehresmann's
theorem, the fibers of $\pi$ are diffeomorphic. Then
$\pi$ is a locally trivial fibration, hence for any fixed $k$ the
$k$-th cohomology of its fibers form a local system, called
{\bf Gauss-Manin local system}.
By the Riemann-Hilbert correspondence, the category of local systems
is equivalent to the category of vector bundles equipped with
a flat connection.

The bundle associated with the  Gauss-Manin local system
is called {\bf the Gauss-Manin bundle}, and the connection
{\bf the Gauss-Manin connection}. 
It can be constructed as follows.

Let $T_\vert M\subset TM$ be the bundle of fiberwise tangent vectors.
By definition, an Ehresmann connection $e$ on $M$ 
is a decomposition $TM= T_\vert M \oplus T_\hor M$,
that is, a choice of such a splitting.
Identifying $T_\hor M$ and the pullback $\pi^*TB$,
we may consider the pullback of a vector field
$X\in TB$ as a vector field $X_e\in T_\hor M$.

A section of a vector bundle associated
with the fiberwise cohomology of $M$
is given by a fiberwise closed differential form 
$\eta\in\Lambda^k M$. The Lie derivative 
$\Lie_{X_e}\eta$
is closed on fibers of $\pi$ for any vector field
on $B$ lifted to a horizontal vector field $X_e$ on $M$.
Indeed, the corresponding diffeomorphisms
map fibers to fibers and fiberwise closed
forms to fiberwise closed forms.

Since different choices $e, e'$ of the
Ehressmann connection result in the vector fields
$X_e, X_{e'}$ which satisfy $Y:=X_e-X_{e'}\in T_\vert M$,
and the form $\Lie_{Y}\eta$ is fiberwise exact,
the cohomology class of the restriction of
$\Lie_{X_e}\eta$ is independent from the
choice of the Ehresmann connection.

Let now $[\eta]$ be the collection of the
cohomology classes of $\eta$ on all fibers of
$\pi$, considered as a section of the Gauss-Manin
bundle. Define 
\begin{equation}\label{_GM_conne_Equation_}
\nabla_X[\eta]:= [\Lie_{X_e}\eta],
\end{equation}
where $[\Lie_{X_e}\eta]$ is the collection
of the cohomology classes of $\Lie_{X_e}\eta$
on all fibers of $\pi$. This formula
defines the Gauss-Manin connection $\nabla$.

\hfill

We needed this observation to prove the following
lemma.

\hfill

\lemma\label{_Gauss_Manin_d_1_vanishes_Lemma_}
Let $\pi:\; M \arrow B$ be a smooth fibration, and
$\eta$ a $p$-form which is closed on the fibers of $\pi$.
Using the pullback map, we consider $\pi^*(\Lambda^*B)$ 
as a subspace in $\Lambda^*(M)$.
Assume that $d\eta$ belongs to 
$\pi^*(\Lambda^2B)\bigwedge \Lambda^{p-1}M\subset \Lambda^{p+1}(M)$.
Then the section of the Gauss-Manin bundle, corresponding
to $[\eta]$, is parallel.

\hfill

\proof
We use the formula \eqref{_GM_conne_Equation_}: 
$\nabla_X[\eta]= [\Lie_{X_e}\eta]$. Cartan's formula gives 
$\Lie_{X_e}\eta= i_{X_e} d\eta + d(i_{X_e}\eta)$, the second term
on the right hand side is exact, and the first vanishes on fibers
because $d\eta\in \pi^*(\Lambda^2B)\bigwedge \Lambda^{p-1}M$.
\endproof

\subsection{Elliptic fibrations on  non-K\"ahler surfaces}

For a version of the following theorem (with a different proof),
see \cite[Theorem 3.17]{_Brinzanescu:bundles_} and \cite{_Brinzanescu:manuscripta_}.

\hfill

\theorem \label{_elli_bundle_Theorem_} 
Let $M$  be a non-K\"ahler compact complex surface
admitting a non-constant 1-dimensional 
family of divisors. Then $M$ 
admits a holomorphic, surjective map $\pi:\; M \arrow S$ 
to a curve, and the general fibers of 
$\pi$ are isomorphic elliptic curves, homologous to 0.
If, moreover, $M$ is minimal, then all fibers of
$\pi$ are elliptic curves.

\hfill

\pstep We may assume that all the divisors in our family are
irreducible curves (if all components in a divisor do not move, the
divisor itself does not move). Let $S$ be the base of the
family. This is a reduced irreducible complex space of dimension
$1$, and the curves in the family with multiplicity $1$ form an
analytic family of $1$-cycles on $M$. By the universal property of
the Barlet spaces, the family is induced from the universal one via
a map $S \to {\goth B}_1(M)$. Replace $S$ with the irreducible
component of ${\goth B}_1(M)$ that contains the image of the
map. Then as shown in \ref{_Barlet_on_surf_Corollary_}, each
connected component of ${\mathfrak B}_1(M)$ is compact. Therefore
$S$ is compact, the support $S^m \subset M \times S$ of our family
is irreducible, and its image $\Psi(S^m) \subset M$ under the
projection $\Psi:M \times S \to M$ is an irreducible closed
subvariety. Since the family is not constant, $\Psi(S^m)$ cannot be
a curve, so $\Psi(S^m)=M$.

All complex curves $C\subset M$ with the fundamental class $[C]\neq
0$ have negative self-intersection by
\ref{_neg_def_h^11_Proposition_}. Therefore if we let $U \subset S$
be the dense open subset parametrizing cycles that are irreducible
with multiplicity one, and let $U^m \subset S^m$ be it preimage in
$S^m$, then the proper map $\Psi:S^m \to M$ is injective on
$U^m$. First of all, this means that $S$ is of dimension exactly
$1$, and $S^m$ is of dimension $2$. Secondly, $\Psi$ is one-to-one
over the complement $M \setminus \Psi(S^m \setminus U^m)$, and since
$S^m \setminus U^m$ is a finite disjoint union of ``bad'' fibers of
the projection $S^m \to S$, $\Psi(S^m \setminus U^m) \subset M$ is
at most one-dimensional. Lastly, note that while these ``bad''
fibers are different as cycles, they might have common irreducible
components, so they might intersect. However, if a curve $C \subset
S^m$ is contracted by $\Psi$, then it must lie in $S^m \setminus
U^m$, so each of its connected components lies in a bad fiber. But
by definition, $\Psi$ is injective on each of the fibers including
bad ones, so $C$ cannot exist.

Thus $\Psi:S^m \to M$ is a proper map with finite fibers that is
one-to-one over the complement to a subvariety of positive
codimension. Since $M$ is smooth, thus normal, $\Psi$ is an
isomorphism, so we obtain a holomorphic map $\pi:M \cong S^m \to S$
to an irreducible compact curve. Taking its Stein factorization, we
may assume that $S$ is normal, thus smooth.

\hfill

{\bf Step 2:}
Let now $C$ be a general fiber of $\pi$.
Using the standard isomorphism
$\pi^*(\Omega^1S)= N^*_\pi M$ between the pullback of the
cotangent sheaf and the conormal bundle to the fibers of
$\pi$ we obtain an exact sequence, sometimes 
called ``adjunction formula'':
\begin{equation}\label{_adjunction_Equation_}
0 \arrow \pi^*(\Omega^1S)\arrow \Omega^1 M \arrow
\Omega^1_\pi M \arrow 0,
\end{equation}
where $\Omega^1_\pi M$ is the bundle of 
holomorphic differentials on fibers of $\pi$. 
This sequence makes sense on the smooth locus
of $\pi$, but we need it only in a neighbourhood
of the generic fiber $C$.

The exact sequence \eqref{_adjunction_Equation_}
gives $K_M\restrict C= K_C$, that is, 
the canonical bundle to $M$ restricted
to $C$ gives the canonical bundle to $C$.
Since the self-intersection of $C$ is zero,
$C$ is homologous to 0 (\ref{_neg_def_h^11_Proposition_}).
This gives $\int_C c_1(K_M)=0$,
hence the degree of $K_C$ is equal 0.
This implies that $C$ is an elliptic curve.

\hfill

{\bf Step 3:}
Return to the map $\pi:\; M \arrow S$ constructed above.
At this step we are going to prove that all its smooth fibers are isomorphic.

Denote by $S_0\subset S$ its smooth locus,
and let $H$ be the Gauss-Manin bundle
associated with the first cohomology of the fibers of $\pi$.
The closed form $\theta$ gives a non-zero
cohomology class in each of the elliptic
curves $\pi^{-1}(s)$ (Step 3), hence it 
gives a section of $H$. Being closed, it is
constant with respect to the Gauss-Manin connection.
The form $d^c\theta= \Theta$ belongs to $\pi^*\Lambda^2 B$.
Then $I(\theta)$ is closed on the fibers of $\pi$. By \ref{_Gauss_Manin_d_1_vanishes_Lemma_},
$I(\theta)$, and hence the Hodge components of $\theta$,
define a parallel section of the Gauss-Manin 
bundle.\footnote{The same conclusion is implied 
by Schmid's fixed part theorem
(\cite[Theorem 7.22]{_Schmid:singularities_}):
for any polarizable variation of Hodge structures over
a quasiprojective base,
the (1,0) and (0,1)-parts of a parallel section $\theta$ are also
constant with respect to the Gauss-Manin connection.}

We obtain a basis $\theta^{1,0}$ and $\theta^{0,1}$
for $H$. Therefore, the variation of Hodge structures 
induced by periods of the elliptic curves on
$H$ is trivial. This implies that the
corresponding elliptic fibration is isotrivial.

Note that isotriviality implies local triviality
by a theorem of Grauert and Fischer
(\cite{_Grauert_Fischer_}), hence
the fibration $\pi$ is locally trivial over $S_0$.

\hfill

{\bf Step 4:} It remains to show that each fiber of $\pi$ contains
an elliptic curve, and nothing else if $M$ is minimal.

Take a critical value $s \in S \setminus S_0$ of the map $\pi$, and
choose a small neighborhood $U \subset s$ that contains no other
critical values, and such that $M_U = \pi^{-1}(U)$ admits a
retraction onto the special fiber $C_s = \pi^{-1}(s)$
(this retraction is constructed in 
\cite{_Persson:degene_}, \cite{_Clemens:degene_},
see also  \cite{_Morrison:Clemens_Schmid_}). Choose also a
point $u$ in the punctured disc $U^o = U \setminus \{s\}$, a point
$m \in C_u$ in its preimage $C_u = \pi^{-1}(u) \subset M_U$, and
denote $M_U^o = \pi^{-1}(U^o) = M_U \setminus C_s \subset M_U$.

Since $C_s \subset M$ is of real codimension at least $2$, the map
$\pi_1(M_U^o,m) \to \pi_1(M_U,m)$ is surjective, and then so is the
map $H_1(M_U^o,\Z) \to H_1(M_U,\Z) \cong H_1(C_s,\Z)$ and its
rational version
\begin{equation}\label{spec.eq}
H_1(M_U^o,\Q) \to H_1(M_U,\Q) \cong H_1(C_s,\Q).
\end{equation}
By the Leray spectral sequence, we have $H_1(M_U^o,\Q) \cong
H_1(C_u,\Q) \oplus \Q\langle \alpha \rangle$, where $\alpha$ is any
class with non-trivial image $\pi(\alpha) \in \Q \cong
H_1(U^o,\Q)$. In particular, it induces a canonical specialization
map $H_1(C_u,\Q) \to H_1(C_s,\Q)$ and its counterpart in homology
with coefficients in $\Z$. We can then take a multisection $\sigma:\tilde{U} \to M_U$
of the projection $\pi:M_U \to U$ --- that is, a finite cover
$\tilde{U} \to U$ of our disc $U$ ramified at $s \in U$ that factors
through $\pi$ --- and let $\alpha = \sigma(1)$ be the image of the
generator $1 \in \Q$ under the map $\sigma:\Q = H_1(\tilde{U}^o,\Q)
\to H_1(M_U^o,\Q)$. Then the map \eqref{spec.eq} annihilates
$\alpha$, so that the specialization map $H_1(C_u,\Q) \to H_1(C_s,\Q)$ is
still surjective.

Next, note that since $U$ is a disc, the K\"ahler form $\omega_S$ is
$dd^c$-exact on $U$, so that $\omega_S=\bar\6\6 \phi$,
for some $\phi\in C^\infty U$. Replacing $\theta^{1,0}$ by
$\gamma:=\theta^{1,0}- \6\phi$, we obtain 
a (1,0)-form $\gamma$ with the same restriction to each fiber, but
now it satisfies 
\[ d\gamma= d \theta^{1,0}- d\6\phi= 
\omega_S-\bar\6\6 \phi=0.
\]
Therefore, $\gamma$ is holomorphic and closed on $M_U$, and as such, it has a
de Rham cohomology class $[\gamma]$.

Denote $P = H_1(C_u,\Z)$, $P' =H_1(C_s,\Z)$, and consider the maps
\begin{equation}\label{per.eq}
P \to P' \to \C,
\end{equation}
where the first arrow is the specialization map, and the second one
is obtained by evaluating $[\gamma]$. We know from Step 2 that $C_u$
is an elliptic curve $E$, so that $P \cong \Z^2$ is a
rank-$2$-lattice, and the map \eqref{per.eq} is the period map that
identifies $E \cong \C/P$. In particular, $P \to \C$ is
injective. Then so is the specialization map $P \to P'$, and since
$P \otimes \Q \to P' \otimes \Q$ is surjective, we see that $P'
\otimes \Q \cong P \otimes \Q$, and $P \subset P'$ is a subgroup of
some finite index $n$.

Denote $G = P'/P$, note that being a quotient of $H_1(M_U,\Z) \cong
P'$, it is also a quotient of $\pi_1$, and let $\eta:M'_U \to M_U$
be the corresponding unramified $n$-fold Galois cover with the
Galois group $G$. By definition, the map $H_1(M'_U,\Z) \to
H_1(M_U,\Z)=P'$ factors through $P$. Consider the map $\pi
\circ \eta:M'_U \to U$, and take its Stein factorization $\pi \circ
\eta = \nu \circ \eta'$ into a finite ramified covering $\nu:U' \to
U$ and a holomorphic map $\pi':M'_U \to U'$ with connected
fibers. Then for any $u \in U^o$, the map $\pi(C_u) \to \pi(M_U) \to
G$ is trivial by construction, so that $\eta$ splits over $C_u$, and
$C'_u = \eta^{-1}(C_u)$ is the disjoint union of $n$ copies of $E
\cong C_u$ transitively permuted by $G$. By definition, these
components correspond to points in the set $\nu^{-1}(u)$, so that
$\nu:U' \to U$ is an $n$-fold unramified cover over $U^o \subset U$.
Since $M'_U$ is
connected, so is $U'$, and then being normal hence smooth, it must
be a disc, and $\nu:U' \to U$ must be the standard $n$-fold cover $z
\mapsto z^n$. Therefore $G \cong \Z/n\Z$ is a cyclic group whose
generator acts on $U'$ via the root of unity.

Now let $C_s'=\eta^{-1}(C_s)$, choose a point $m \in C'_s$, and
consider the Albanese-type map $\Alb:M'_U \to E=\C/P$ sending a
point $m' \in M_U'$ to the integral $\int_l \eta^*\gamma$, where
$l:[0,1] \to M'_U$ is any path connecting $m$ and $m'$ (modulo $P$,
the integral does not depend on the choice of the path). Then $\Alb$
identifies each component of the curve $C'_u$, $u \in U^o$ with $E$,
and the product map $\Alb \times \pi':M'_U \to E \times U'$ is a
proper map that is an isomorphism over $\nu^{-1}(U^o)$. Thus it is a
surjective rational map, $\pi_1(M'_U,m) \cong P$ is abelian, and
$H^1(E,\Z) \cong H^1(E \times U',\Z) \to H^1(M'_U,\Z) \cong
H^1(C'_s,\Z)$ is an isomorphism. Moreover, $C'_s \to E$ is a proper
surjective map, so that $C'_s$ has a component $C$ that dominates
$E$, thus has genus $\geq 2$, and possibly some other
components. But if we let $C''_s$ be the normalization of $\C'_s$,
then the map $H^1(E,Z) \cong H^1(C'_s,\Z) \to H^1(C''_s,\Z)$ is
surjective. Therefore $H^1(E,\Z) \cong H^1(C,\Z)$ so that $C \cong E$,
while all the other components are rational curves with no $H^1$. If
$M$ is minimal, the other components do not exist and $M'_U \cong E
\times U'$.

To finish the proof, it remains to notice that $G \in \Z/n\Z$ acts
on $M'_U$ without fixed points, and since it fixes the only
non-rational component $C \cong E$ of $C'_s$, it must also act on
$E$ without fixed points. Therefore the generator $1 \in P'/P \cong
\Z/n\Z$ acts by shifts by an $n$-torsion element $e \in E = \C/P$,
the map $P' \to \C$ is injective, so that $P' \cong \Z^2$, and $C_s
= C'_s/G$ consists of an elliptic curve $E' = E/G = \C/P'$ and some
rational curves (or nothing in the minimal case).
\endproof
 
\hfill

\corollary\label{_fibra_Corollary_}
Let $M$ be a non-K\"ahler, compact, complex 
surface with $b_1(M)>3$. Then $M$ admits an elliptic fibration.

\hfill

\proof
 \ref{_H^1_odd_Theorem_} implies that $\dim H^{1,0}(M)= \frac{b_1(M)-1}2$, and hence
$\dim H^0(\Omega^1(M))\geq 2$.
By \ref{_prod_hol_1_forms_Proposition_},
all globally defined holomorphic 1-forms on $M$ are pointwise proportional.
Consider the rank 1 sheaf $L\subset \Omega^1(M)$ generated by
holomorphic 1-forms. If $\dim H^{1,0}(M)=\dim H^0(L) > 1$, the zero
divisors of the sections of this sheaf form a continuous family
of curves on $M$, and \ref{_elli_bundle_Theorem_}
can be applied.
\endproof

\hfill

\remark\label{kodaira.rem}
Note that the result we have proven does not imply that the map
$\pi$ is a submersion. Indeed, $\pi$ may have multiple fibers. The
basic local example is given by Kodaira log-transform. Take an
elliptic curve $E = \C/P$, and a torsion point $e \in E$ of some
order $n$, and consider the quotient $M_U(E,e) = (E \times
U')/(\Z/n\Z)$, where $U'$ is the unit disc, and the generator $1 \in
\Z/n\Z$ acts by translation by $e$ on $E$, and by multiplication with the
$n$-th root of unity on $U'$. Then the action has no fixed points,
so that $M_U(E,e)$ is smooth, and we have a holomorphic map
$M_U(E,e) \to U = U'/(\Z/n\Z)$, where $U$ is again a disc. As we saw
in the proof of Step 4 of \ref{_elli_bundle_Theorem_},
this local example is universal: whenever $M$ is minimal,
all the varieties $M_U = M'_U/(\Z/n\Z)$ are of this type.

\subsection{Isotrivial elliptic fibrations}

\proposition\label{_isotrivial_Proposition_}
Let $M$ be  a non-K\"ahler compact complex 
surface, admitting an elliptic fibration $\pi:\; M \arrow S$. 
Denote by $C$ a general fiber of
$M$. Then there is a natural holomorphic action
of $C$ on $M$, transitive on fibers and free on 
non-multiple fibers.

\hfill

\proof
Two pairs of points $(a, b)$ and $(a_1, b_1)$ on
a smooth genus 1 curve $T$ are called {\bf rationally equivalent}
if there exists a translation automorphism of $T$ mapping
$(a, b)$ to $(a_1, b_1)$. Clearly, 
the set $\Jac(T)$ of such pairs up to rational equivalence is a curve 
isomorphic to $T$, with a fixed point $(x,x)$,
which we can consider as the group unit.
We can consider $T$ as a torsor over the group $\Jac(T)$. \footnote{A
torsor over a group $G$ is a set with a free and transitive $G$-action.}

Consider the relative Jacobian $\Jac(M)$ of $M$ over $B$
(\cite{_Lopes_Martin:relJac_}), that is,
a quotient of the product $M\times_B M$ by rational
equivalence: two points $x, y\in M\times_B M$
are rationally equivalent if they can be
connected by rational curves.  Locally, we can interpret 
the relative Jacobian using the relative Albanese map
(see  \ref{kodaira.rem} and Step 4 of \ref{_elli_bundle_Theorem_}).

The map $\pi:\; \Jac(M)\arrow S$
is a locally trivial fibration with a globally 
defined section, hence it is trivial. The natural
fiberwise action of $\Jac(M)=C\times S$ on $M$ defines 
the action on $M$, which is free and transitive on 
general fibers. This action is transitive on all
fibers. To see this, consider the set
\[
\{(x, y, c) \in M\times_S M\times_S \Jac(M) ,  \ \ |\ \  (x,y)\sim c\}
\]
This set is complex analytic in $M\times_S M\times_S \Jac(M)$.
By Remmert's proper mapping theorem, its projection $V$ to $M\times_S M$ is
closed. Since $V$ contains all regular 
points of the projection $M\times_S M\arrow S$.
one has $V=M\times_S M$.
\endproof

\hfill

\definition
Let $M$ be a complex manifold equipped with
an action of a compact complex torus $T$, 
with the orbits of the same dimension. Assume that
 the quotient map $\pi:\; M \arrow S$ is well defined. Then 
$\pi:\; M \arrow S$ is called {\bf an isotrivial
toric fibration}, and {\bf an isotrivial
elliptic fibration} if $\dim_\C T =1$.

\hfill

\remark
In some literature, ``isotrivial elliptic fibration''
is a fibration with all smooth fibers isomorphic
to the same elliptic curve and the only singularities in 
multiple fibers. In \cite{_Brinzanescu:bundles_}
and \cite{_Brinzanescu:manuscripta_}, the same notion is 
called ``a quasi-bundle''.

\hfill

\remark\label{_orbifold_Remark_}
Let $\pi:\; M \arrow S$ 
be an isotrivial elliptic fibration over a curve $S$, with
fiber $C$. We associate with $\pi$ the following orbifold 
structure on $S$. For any multiple fiber $R=\pi^{-1}(s)$ of $\pi$, consider
the group $\Gamma_R\subset C$, obtained as the kernel of the
natural action of $C$ on $R$. Taking the quotient
$M/\Gamma_R$, we obtain another fibration, which is
locally trivial in a neighbourhood of $s$. A smooth section
of this fibration in a neighbourhood of $s$ gives 
a $\Gamma_R$-invariant multisection $\tilde U \arrow U$ 
of $\pi:\; \pi^{-1}(U) \arrow U$. Then $U$ is obtained
as a finite quotient of $\tilde U$, and the lifting
of $\pi$ to $\tilde U$ is locally trivial.

\hfill

\remark
We obtained that any isotrivial 
elliptic bundle $\pi:\; M \arrow S$  over a curve 
defines an orbifold structure on this curve,
and $\pi$ is locally trivial in this orbifold structure.
The same argument would work for any base $S$
regardless of dimension.

\hfill


To simplify the notation and the arguments,
we shall deal with the smooth orbifold fibrations in the same
way as for the smooth ones, and use the standard
terminology instead of adding ``orbifold smooth''
everywhere. Most of the standard results
and constructions in smooth category are extended 
to the orbifold category in the usual way; for the only exception of
importance to us, see below \ref{_c_1_isotri_Remark_}.

\hfill

The topology and geometry of isotrivial toric fibrations
and principal toric fibrations was explored in some
depth in \cite{_Hofer:remarks_}.
For our present purposes we shall need the following
theorem. Notice that the isotrivial toric
fibration can be considered as a locally trivial principal
toric fibration if we work in the orbifold category.
Following \cite{_Hofer:remarks_}, we define the
Chern classes associated with isotrivial toric
fibrations as follows.

\hfill

\theorem\label{_c_1_toric_Theorem_}
Let $T(S)$ be the sheaf of $T$-valued smooth functions on
$S$. Let $\tilde T=\R^n$ the universal cover of $T$, and
denote by $\tilde T(S)$ the sheaf of smooth $\tilde T$-valued
functions on $S$. Denote by $\Z^n(S)$ the constant sheaf
with fiber $\Z^n$.
Then the following exact sequence of sheaves
\[
0 \arrow \Z^n(S) \arrow \tilde T(S)\arrow  T(S)\arrow 0
\]
gives an exact sequence in cohomology
\[
0 = H^1(S,\tilde T(S)) \arrow H^1(S,TS) \arrow H^2(S,\Z^n(S))
\arrow H^1(S,\tilde T(S)) =0
\]
and a bijective correspondence between
the set $H^1(T(S))$ of principal toric fibrations 
and the cohomology group $H^2(\Z^n(S))= H^2(S, \Z)^n$.

\proof
Clear (see \cite{_Hofer:remarks_} for details). \endproof

\hfill

\definition
Let  $\pi:\; M \arrow S$ be an isotrivial toric fibration,
with $T=\R^n/\Z^n$ its fiber. The $n$ cohomology classes
associated to $\pi$ as in \ref{_c_1_toric_Theorem_}
are called {\bf the Chern classes of the isotrivial toric
  fibration}.

\hfill

\remark \label{_c_1_isotri_Remark_}
In the case when $\pi:\; M \arrow S$ 
is an isotrivial elliptic fibration over a 
compact curve, one has $H^2(S, \Z)=\Z$. Therefore
a toric fibration is uniquely determined by the
vector $c_1(\pi) \in \Z^2= H^2(S, \Z)^2$. Recall that a vector
$v$ in a lattice $\Z^n$ is called {\bf primitive}
if $v$ is not divisible by an integer $n>1$.
For any primitive
vector  $v$ in a lattice  $\Lambda=\Z^2$, 
there exists $w\in \Lambda$ such that 
$\Lambda = \langle w, v\rangle$. Therefore,
the Chern classes of an 
isotrivial elliptic fibration over a curve
are determined (up to an automorphism of $\Lambda$) 
by the biggest $n \in \Z$ such that $c_1(\pi)$ is
divisible by $n$. Using this observation, we shall understand
$c_1(\pi)$ as a non-negative integer; it is equal to zero when
$v=0$ and equal to the largest integer divisor of $v$ when $v\neq
0$. Note that if $S$ is an orbifold, then the first Chern class is
still perfectly well-defined and classifies isotrivial elliptic
fibrations, but it is no longer necessarily true
that $H^2(S,\Z) = \Z$ --- it might also have some torsion.

\hfill

Now we shall prove that an isotrivial elliptic fibration 
$\pi:\; M \arrow S$ over a curve with non-torsion $c_1(\pi)$ 
is always obtained as a quotient of the total space
of an ample $\C^*$-bundle. We need the following preliminary
topological statement. 

\hfill

\claim\label{_c_1_via_cover_Equation_}
Let $\pi:\; M \arrow S$ 
be an isotrivial elliptic fibration over a compact curve $S$,
and $C$ its general fiber. Assume that $c_1(\pi)$ is not torsion.
Then the rank of the natural map $H_1(C) \arrow
H_1(M)$ is 1. Moreover, there exists
a normal subgroup $G\subset \pi_1(M)$ with $\pi_1(M)/G=\Z$ such that
the corresponding $\Z$-cover is infinite on $C$.

\hfill

\proof
Consider the exact sequence
\begin{equation}\label{_Leray_exact_Equation_}
H_2(S,\Q)\stackrel \delta \arrow  H_1(C,\Q) \stackrel \psi \arrow H_1(M,\Q) \arrow
H_1(S,\Q)\arrow 0
\end{equation}
obtained from the Leray spectral sequence of the fibration.
It is easy to see that $\delta$ is dual to the Chern
class of $\pi$. Therefore, the map $\psi:\; H_1(C,\Q) \arrow
H_1(M,\Q)$ has rank 1 when $c_1$ it not torsion and 2 when
it is. Therefore, there exists an element
$v\in \pi_1(C)$ such that its image in $\pi_1(M)$ 
has infinite order.  Take a homomorphism
$H_1(M,\Z) \arrow \Z$ which is non-zero on $\psi(H_1(C))$.
Since $H_1(M)= \frac{\pi_1(M)}{[\pi_1(M),\pi_1(M)]}$,
the map $v$ defines a group homomorphism
$\pi_1(M) \arrow \Z$ which is non-trivial on
$\pi_1(C)= H_1(C,\Z)$.
\endproof

\hfill

\proposition\label{_elli_via_line_bundle_Proposition_}
Let $\pi:\; M \arrow S$ 
be an isotrivial elliptic fibration over a compact curve $S$,
and $C$ its general fiber. Consider a $\Z$-cover $\phi:\;\tilde M
\arrow M$ such that $\phi^{-1}(C)\arrow C$ is an infinite cover.
Then $\tilde M$ is the total space of a principal $\C^*$-bundle
associated with a line bundle $L$ on $M$. Moreover, the surface $M$ is obtained
as the quotient of the total space $\Tot(L^\circ)$ of non-zero
vectors in $L$ by a holomorphic automorphism  $q:\; S \arrow S$
acting on $L$ equivariantly and linearly on all fibers.
Finally, the automorphism $q:\; S \arrow S$
has finite order, and $|c_1(L)|$ is equal to $c_1(\pi)$.

\hfill

\proof
Since $\phi^{-1}(C)\arrow C$ is a $\Z$-cover, the 
space $\tilde M$ is a $\C^*$-bundle over $S$. Denote the
associated vector bundle by $L$. Then $\tilde M=\Tot(L^*)$,
and $M$ is obtained from $\Tot(L^\circ)$ as a quotient
by an cyclic group of automorphisms generated by
an automorphism $\tilde q:\; \tilde M \arrow \tilde M$
commuting with the $\C^*$-action. Then $\tilde q$ defines
a holomorphic automorphism of $S$.
However, all automorphisms of $\C^*$ commuting with the 
$\C^*$-action are given by a multiplication by a number. 
Therefore, $\tilde q$ gives a holomorphic section 
$\lambda\in \Hom(L, q^*(L))$. The fibers of
$\pi:\; M \arrow S$ are compact, hence $q$ has to be of finite order.

The Chern class of 
$L$ is obtained from the $U(1)$-bundle associated with
$L$, and this bundle is homotopy equivalent to $\tilde M$,
hence $c_1(L)$ is equal to $c_1$ of the circle bundle
associated with the $\C^*$-bundle $\tilde M \arrow S$.
However, $c_1$ of $\tilde M \arrow S$ is by construction
equal with $\pm c_1(\pi)$ (the $\pm$ sign is due to the
ambiguity of the definition of $c_1$ of the toric 
fibration, see \ref{_c_1_isotri_Remark_}).
\endproof

\subsection{Blanchard theorem and LCK structure on 
non-K\"ahler elliptic surfaces}

The main result of the present section
is the following theorem, proven at the end of this subsection.

\hfill

\theorem\label{_elli_then_Vaisman_Theorem_}
(\cite[Theorem 1]{_Belgun_})
Let $M$ be a compact, non-K\"ahler, minimal surface
admitting an elliptic fibration. Then $M$ is 
Vaisman.

\hfill

We use Blanchard theorem (see \cite{_Blanchard_} or 
\cite{_Rogov:Iwasawa_}), which is applied to the present case as follows.

\hfill

\theorem\label{_Blanchard_Theorem_}
Let $\pi:\; M \arrow S$ be an isotrivial 
elliptic fibration. Then $M$ is K\"ahler if and only if
$c_1(\pi)$ is torsion.

\hfill

\proof
Suppose that $c_1(\pi)$ is torsion, and consider the exact sequence
\[
0 \arrow H^1(S,\Q) \arrow H^1(M,\Q) \arrow H^1(C,\Q) \stackrel
{\delta^*} \arrow H^2(S,\Q)
\]
obtained from the Leray spectral sequence.
As shown in the proof of \ref{_c_1_via_cover_Equation_},
one has $\delta^*=0$ if and only if $c_1(\pi)$ is torsion.
However, if $\delta^*\neq 0$, the pullback
$\pi^*(\omega_S)$ of the volume form is exact,
which is impossible when $M$ is K\"ahler by 
\ref{_exact_pos_non-Kahler_Remark_}. 

Assume now that $c_1(\pi)$ is torsion, or equivalently, that
$\delta^*=0$. Using \ref{_elli_via_line_bundle_Proposition_},
we obtain $M$ as a quotient of a total bundle
$\Tot(L^\circ)$ by an automorphism acting on the base $S$ as
$\phi:\; S \arrow S$ and on
the fibers as a constant endomorphism
$\lambda:\; L \arrow \phi^* L$, $|\lambda|\neq 0$. Since $c_1(L)=
\pm c_1(\pi)$ is torsion,
this bundle admits a flat Hermitian connection.
Denote by $\rho$ its monodromy action.
Since $\rho$ acts on the fibers of $L$ isometrically,
we can chose a $\rho$-invariant K\"ahler metric $\omega_F$
on  the fibers of $L^\circ/\lambda$, which are elliptic curves.
Extending $\omega_F$ to $M$ using the flat Ehresmann
connection on $\pi:\; M \arrow S$, we obtain a 
closed, positive (1,1)-form $\hat \omega_F$ on $M$ which is
strictly positive on the fibers of $\pi$;
the sum $\omega:=\hat\omega_F+ \pi^*(\omega_S)$ is
a positive, closed, strictly positive almost
everywhere (1,1)-form. Then $\omega^2>0$,
in contradiction with the negativity of the intersection
form on $H^{1,1}(M)$ (\ref{_neg_def_h^11_Proposition_}). 
\endproof

\hfill

Let us prove \ref{_elli_then_Vaisman_Theorem_} (see also \cite[Theorem 3.5 and the subsequent Remark]{_Vaisman:gen_hopf_} for the regular case, and also \cite{_Vuli:arxiv_}).
If $M$ is non-K\"ahler, it is obtained as a quotient
of $\Tot(L^\circ)$  by $\lambda:\; L \arrow \phi^* L$,
where $c_1(L)\neq 0$ by \ref{_Blanchard_Theorem_}. 
Replacing $L$ by its dual
bundle if necessary, we may assume that $L$ is ample.
Choose a  Hermitian structure on $L$ such that
its curvature is a positive (1,1)-form on $S$,
and let $\psi\in C^\infty(\Tot(L^\circ))$ be the  
function $\psi(v) = |v|^2$. 

We want the isomorphism $\lambda:\; L \arrow \phi^* L$
to have constant length. Since $\phi$ is of finite order,
we can always represent the curvature $\Theta$ of $L$ by
a $\phi$-invariant form. Using the $dd^c$-lemma as usual,
we choose a metric on $L$ with curvature $\Theta$.
Then the curvature of $\Hom(L, \phi^*L)$ is $\phi^*(\Theta)-\Theta=0$.
 Therefore, $\Hom(L, \phi^*L)$ is a flat unitary bundle; any section
$\lambda\in \Hom(L, \phi^*L)$  gives a holomorphic
section of a flat unitary bundle, hence it
is of constant length.

Then $dd^c(\psi)$
is a K\"ahler form on $\Tot(L^\circ)$
(\cite[(15.19)]{_Besse:Einst_Manifo_}).
Then $\lambda$ acts on $\Tot(L^\circ)$ by holomorphic
homotheties, hence the quotient $M = \Tot(L^\circ)/\langle \lambda\rangle$
is LCK. Also, this quotient is Vaisman, because the
standard action of $\C^*$ on $\Tot(L^*)= \tilde M$
is by holomorphic homotheties, \cite{_Kamishima_Ornea_}.


\section{Class VII surfaces}


In this section we prove the following theorem.

\hfill

\theorem\label{_elli_or_class_VII_Theorem_}
Let $M$ be a non-K\"ahler compact minimal complex surface. 
Then $M$ is isotrivial elliptic or belongs to class VII
(or both). 

\hfill

\pstep 
As $b_1(M)>0$, and $H_1(M, \Z)= \pi_1(M)/[\pi_1(M),  \pi_1(M)]$,
there exists a subgroup of any given finite  index $r$ in $\pi_1(M)$.
Denote the corresponding cover by $\sigma:\; M_1\arrow M$.
Then $M_1$ is a compact manifold. 
Since $\chi(\calo_{M})= \frac{c_1^2+c_2}{12}$ is expressed through the
curvature, one has $\chi(\calo_{M_1})= d \chi(\calo_{M})\geq d$
for any $d$-sheeted cover $\sigma:\; M_1\arrow M$.
Unless $\chi(\calo_{M})=0$, we can find a cover $M_1$ with
$\chi(\calo_{M_1}) < -3$ 
or $\chi(\calo_{M_1}) > 3$. In the first case, $M_1$ is elliptic by
\ref{_fibra_Corollary_}. 
In the second case, the canonical bundle
$K_{M_1}$ satisfies $\dim H^0(K_{M_1})\geq 2$, hence it
has sections which give a continuous family of divisors.
Therefore, for $\chi(\calo_{M})\neq 0$, the finite cover $M_1$ has a 
continuous family of divisors. Then
$M$ also has such a continuous family, and it is elliptic by
\ref{_elli_bundle_Theorem_}.

\hfill

{\bf Step 2:} Assume that $b_1(M)=3$, but $M$ is not elliptic. 
Then $\chi(\calo_{M})=0$, hence $h^{0,1}(M)=2$ and $h^{0,2}(M)=1$,
and the same is true for all non-ramified covers of $M$.
We are going to prove that $M$ contains no curves.

First, we prove that $c_1(M)^2=(K_M)^2=0$. Otherwise $c_1(M)^2 < 0$,
so $K_M$ is non-trivial, and since $h^{0,2}(M)=1$, we have $K_M =
\calo(D)$ for some effective divisor $D=\sum_ia_iD_i$. If for some
$i$ we have $(D_i)^2 < 0$ and $(K_M \cdot D_i) < 0$, then $D$ is a
$(-1)$-curve, and this is not possible since $M$ is minimal. If
$(D_i)^2=0$, then $D_i$ is homologous to $0$ by
\ref{_neg_def_h^11_Proposition_}. Thus in any case, $(K_M \cdot
D_i)=0$, and then $(K_M)^2 = \sum_ia_i(K_M,D_i)=0$, a contradiction.

Note that $c_1(M)^2=0$ implies that
$$
0=\chi(\calo_{M})= \frac{c_1(M)^2+c_2(M)}{12}= \frac{c_2(M)}{12} =
\frac{\chi(M)}{12},
$$
hence $b_2(M)= b_1(M) + b_3(M) - b_0(M) - b_4(M) = 6-2=4$. Moreover,
the same holds for any unramified cover $M_1$ of $M$.

Now assume that $M$ contains an irreducible curve $C$. Then by the
adjunction formula, its arithmetic genus is $p_a(C) =
\frac{(C)^2+(K_M \cdot C)}{2}+1$, and since $(C)^2 \leq 0$ and
$(K_M)^2=0$, we have $(K_M \cdot C)=0$, hence $p_a(C) \leq 1$. Then
$C$ is either smooth rational, or rational with a single node, or
rational with a single cusp, or smooth elliptic, and in any case,
$b_1(C) \leq 2$. Since $b_1(M)=3 > b_1(C)$, we have an unramified
cover $M_1 \to M$ of any degree $d$ that splits over $C$, thus
contains non-intersecting curves $C_1,\dots,C_d$ with the same
self-intersection $l=(C_i)^1$, $i=1,\dots,d$.

Then if $l < 0$, the classes of the curves $C_i$ span a
$d$-dimensional subspace in $H^2(M_1,\Q)$, so that $d \leq b_1(M_1)
= 4$, and this is a contradiction since $d$ was arbitrary.

Otherwise $l=0$, and all the curves $C_i$ are homologous to $0$. Let
$D$ be their union, and consider the short exact sequence
\begin{equation}\label{_exact_O_curves_Equation_}
0 \arrow \calo(-D) \arrow \calo_{M_1} \arrow \bigoplus_i \calo_{C_i}\arrow 0.
\end{equation}
The corresponding long exact sequence
\[ \C = H^0(\calo_{M_1}) \arrow \C^d = \bigoplus_i H^0(\calo_{C_i})
\arrow H^1(\calo(-D)).\]
shows that that $\dim H^1(M_1,L) \geq d-1$, where $L:=\calo(-D)$. However,
$\chi(L)=\chi(\calo_{M_1})=0$ because $c_1(L)=[D]=0$, and for any line
bundle $L'$ on $M_1$, we have $\dim H^0(M_1,L') \leq 1$ -- otherwise
$M_1$, hence also $M$ carries a non-trivial family of divisors and is elliptic by
\ref{_elli_bundle_Theorem_}. Thus $\dim H^0(M_1,L)$ and $\dim
H^2(M_1,L) = \dim H^0(M_1,K_{M_1} \otimes L^*)$ are at most $1$, and
$d \leq 3$. This is again a contradiction.

\hfill

{\bf Step 3:}
We can prove now that all surfaces 
with $b_1(M)=3$ and $\chi(\calo_M)=0$ are elliptic.
Arguing by absurd, we may assume $M$ is non-elliptic,
but in this case $M$ has no complex curves (Step 2).
Consider the space $W\subset H^1(M)$ generated by
holomorphic and antiholomorphic forms.
Let $\alpha\in \Lambda^{1,0}(M)$ be a holomorphic 1-form
generating the space $H^{1,0}(M)$. Fix $x\in M$
and define the Albanese map by taking
$y\in M$ to $\int_\gamma \alpha$, where
$\gamma$ is a path connecting $x$ to $y$.
This map depends on the choice of the path $\gamma$,
hence it gives a map $\Alb:\; M \arrow \C /\Lambda$,
where $\Lambda\subset \C$ is the period lattice of $\alpha$,
that is, the set  $\{\int_v\alpha \ \ |\ \ v \in H_1(M,
\Z)\}$.  By 
\ref{_W_rational_Corollary_}, $W$ is rational.
Therefore, the integral $\int_v\alpha$
vanishes on one of the generators in $H_1(M,\Z)$,
and  $\Lambda$ is a rank 2 lattice in $\C$.
It is easy to see that it is discrete. 
Then the Albanese map gives a fibration
$M \arrow \C/\Lambda$, and $M$ contains
a continuous family of curves, hence it is
elliptic.

\hfill

{\bf Step 4:} It remains to show that $M$ is 
of class VII or elliptic if $\chi(\calo_{M_1})=0$.
When $b_1(M)=3$, the surface $M$ is elliptic (Step 3), and the same
holds if $b_1(M) \geq 5$ by \ref{_fibra_Corollary_}. 
Therefore, we may assume that $b_1(M)=h^{0,1}(M)=1$ and
the same is true for all covers $M_1$ of $M$, ramified or
unramified. Unless $M$ is elliptic, we have
$\chi(\calo_M)=0$, hence $h^{0,2}(M)=0$.
To prove that $M$ is class VII one has
to check that $H^0(K_M^n)=0$ for any $n>0$.

If $K^n_M$ admits a non-zero section, then $K_{M_1}$
admits a non-zero section for some
finite ramified cover $\sigma:\; M_1\arrow M$. Indeed, given a 
section $\alpha$ of $K^n_M$, with zero divisor $C$,
locally around a point $x\in M\backslash C$, one can take
its $n$-th degree root $\sqrt[n]{\alpha}$ and obtain a section
of $K_M$. The sheaf of such sections is finite and locally constant
on $M \backslash C$,
hence it becomes trivial after passing to a finite cover 
$\widetilde {M\backslash C} \arrow M\backslash C$. 
In a neighbourhood $U$ of $x\in C$, we interpret
$\sqrt[n]{\alpha}$ as a multivalued function with values in $K_M$,
which can be trivialized in $U$. The graph of $\sqrt[n]{\alpha}\subset U\times C$
can be glued to $\widetilde {M\backslash C}$, giving a 
ramified cover $\sigma:\; M_1\arrow M$ and a section of $K_{M_1}$.
Then $h^{0,2}(M_1)\neq 0$, hence $M_1$ is elliptic (Step 3).
Therefore, $M=\sigma(M_1)$ contains a continuous family
of holomorphic curves; it is elliptic by \ref{_elli_bundle_Theorem_}.
Isotriviality of the elliptic fibration follows from
\ref{_isotrivial_Proposition_}. 
\endproof

\hfill


\section{Brunella's theorem: all Kato surfaces are LCK}
\label{_Brunella_proof_}


The Kato surfaces are also called ``Global spherical shell surfaces''.
By definition a Kato surface is a surface $M$ which possesses a 
{\bf global spherical shell}, that is, an open subset $U\subset M$
which is biholomorphic to a neighbourhood of $S^3$ in $\C^2$
and such that $M \backslash U$ is connected.
In this section we prove Brunella's theorem,
showing that all Kato surfaces are locally
conformally K\"ahler. 

For the original definition and early works on Kato surfaces see
\cite{_Kato:announce_,_Kato:Kinokunya_,_Kato:sugaku_,_Dloussky:Kato_},
and for Brunella's original proof,
\cite{_Brunella:Enoki_} and \cite{_Brunella:Kato_}.

It is not hard to see that all Kato surfaces are 
deformations of a blow-up of a Hopf surface. To see
this, we consider the following explicit 
construction of Kato surfaces.

Let $M$ be a Kato surface, and $S\subset U \subset M$ 
the corresponding 3-sphere. Consider the map $\chi:\; \pi(M) \arrow \Z$
mapping a path $\gamma$ to the intersection index $\gamma \cap S$.
Clearly, $\chi$ is a group homomorphism. Denote the
corresponding $\Z$-cover by $\tilde M$, and the preimages
of $S$ in $\tilde M$ by $S_i$, $i\in \Z$
(these preimages can be enumerated, because
the deck transform group $\Z$ acts on the
set of preimages of $S$ freely and transitively).
Denote by $M_i$ the subset of $\tilde M$
situated between $S_i$ and $S_{i-1}$.
Clearly, each $M_i$ is a fundamental domain of
the deck transform action.

\begin{figure}[ht]
\begin{center}\includegraphics[width=0.45\linewidth]{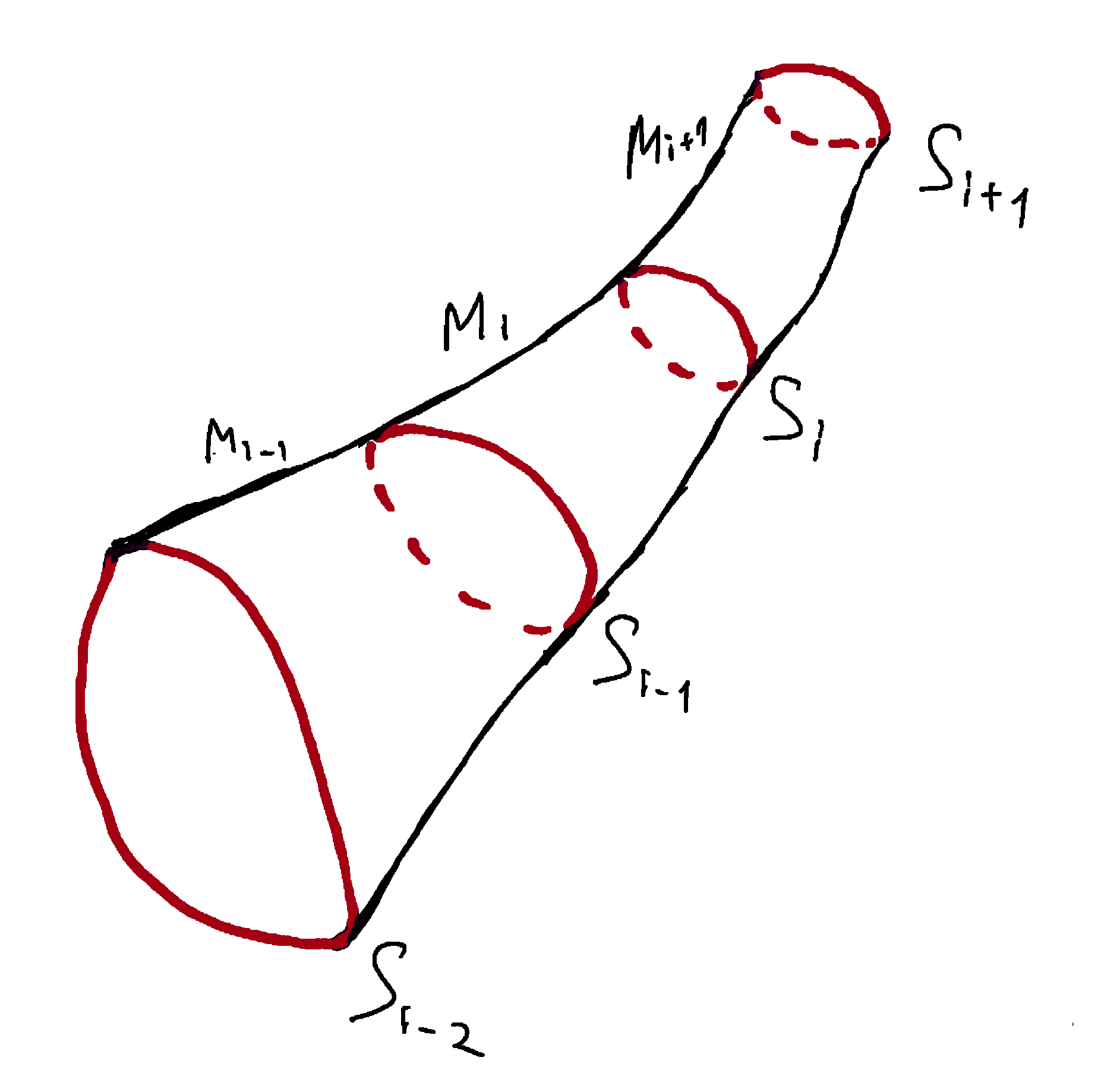}\\
{\em Action of $\Z$ on $\tilde M$, spherical shells and
  fundamental domains}
\end{center}
\end{figure}

Each $M_i$ has two boundary components, with 
$S_i$ pseudo-convex and $S_{i-1}$ pseudo-concave.
Gluing a ball $B$ to $S_{i-1}$, we obtain a manifold
$\hat M_i= B \coprod_{S_{i-1}} M_i$ which is compact and has strictly
pseudoconvex boundary. By Grauert's solution of Levi problem
(\cite[Section 2, Theorem 1]{_Grauert:Levi_}), $\hat M_i$ is holomorphically convex.
Using Remmert's reduction theorem \cite{_Remmert:reduction_}, 
we obtain a proper map with connected fibers
$p:\; \hat M_i \arrow X$, where $X$ is Stein.
Since the neighbourhood of the boundary of
$X$ is biholomorphic to a neighbourhood of $S^3\subset B$,
Hartogs theorem implies that $X$ is in fact biholomorphic
to $B$. Then $p:\; \hat M_i \arrow X$ is a bimeromorphic,
holomorphic map to a ball, and is obtained by a sequence
of a blow-ups.

Following Dloussky \cite{_Dloussky:Kato_}, we define the 
{\bf Kato data} on a closed ball $B\subset \C^2$
as a bimeromorphic, holomorphic map $\hat B \arrow B$,
together with an open subset $B_0\subset \hat B$
and a biholomorphism from $B_0$ to an open ball.
Then the complement $\hat B \backslash B_0$
has two smooth boundary components, which are both isomorphic
to $S^3\subset B$ and can be glued together to obtain
a compact complex surface. 

\hfill

We have just proved the following theorem, originally
due to Ma. Kato. (\cite{_Kato:announce_}). 

\hfill

\theorem
Let $M$ be a Kato surface, and $S\subset M$ its global
spherical shell. Then $M$ can be obtained
from the Kato data by gluing two boundary components
of $\hat B \backslash B_0$ as above.
\endproof

\hfill

\remark\label{_center_Kato_Remark_}
Using the same arguments, it is possible to show
that the Kato surface $M$ is not minimal unless
$\hat B$ is obtained by blowing up a point $x_0\in B$,
then blowing up points which lie on exceptional divisors
(\cite{_Dloussky:private_}).
Since a blow-up of an LCK manifold is LCK, we can
always assume that $M$ is minimal and $\hat B$ is 
obtained from $B$ by successive blow-ups in $0\in B$.

\hfill

\remark\label{_domain_C^2_Remark_}
Suppose that the Kato data satisfy the
assumptions in \ref{_center_Kato_Remark_}.
Then in the construction of the Kato surface $M$ from
the Kato data, we can always replace the open ball $B$ of radius
1 by a ball $B(r)\subset B$ of radius $r$ and
the ball $B_0$ by its image $\Psi(B(r))$.
Clearly, the resulting Kato surface is
biholomorphic to $M$. Similarly, the
ball $B$ can be replaced by any holomorphically
convex domain $U\subset \C^2$ with smooth boundary containing the origin $0$.

\hfill

\definition
Choose $\epsilon >0$, and let 
$\max_\epsilon:\; \R^2\arrow \R$ be a smooth, convex
function, monotonous in both variables,
which satisfies $\max_\epsilon(x, y) = \max(x,y)$
whenever $|x-y|>\epsilon$. Then $\max_\epsilon$ is called
{\bf a regularized maximum} (\cite{_Demailly_1982_}).
It is easy to see that the regularized maximum
of two plurisubharmonic functions is
again plurisubharmonic. This construction
allows one to ``glue'' plurisubharmonic functions
and the corresponding K\"ahler metrics.

\hfill

Now we can prove Brunella's theorem.

\hfill

\theorem\label{_Brunella_Theorem_}
Let $M$ be a Kato surface. Then $M$ is locally conformally
K\"ahler.

\hfill

\pstep Let $\pi:\; \hat B \arrow B$ and $B_0\subset \hat B$ be the
Kato data. Denote by $\Psi:\; B\arrow B_0$ the corresponding
biholomorphic equivalence. Choose a K\"ahler metric $\hat \omega$ on
$\hat B$.

Brunella's theorem is proved by finding a metric
$\hat\omega$ on $\hat B$ with the following automorphic
condition. Consider the space $\tilde M$ obtained by
gluing $\Z$ copies of $\hat B \backslash B_0 = M_i$ as
above. A K\"ahler metric $\tilde \omega$ on $\tilde M$ is called
{\bf $\Z$-automorphic} if the deck transform group
mapping $M_i$ to $M_j$ acts on $(\tilde M, \tilde\omega)$
by homotheties. To obtain such a form we need to find
a K\"ahler form $\hat \omega$ on $\hat B$ 
such that $\hat\omega\restrict{B_0}$ is equal
to $\Psi^*\hat\omega$ in a neighbourhood of
the boundary of $B_0$. If this is true, $\Psi$ acts
by homotheties in a neighbourhood of
$S$. Then the restriction of $\hat\omega$
to $\hat B \backslash B_0 = M_i$
can be extended to a $\Z$-automorphic K\"ahler form on
$\tilde M = \bigcup_{i\in \Z} M_i$.

This is the strategy we follow, except that we replace
$B$ by another strictly pseudoconvex domain as
in \ref{_domain_C^2_Remark_}.

\hfill

{\bf Step 2:}
Using the local $dd^c$-lemma, we can find a smooth
function $\phi$ on $B_0$ such that
$dd^c\phi=\hat\omega\restrict{B_0}$.
Adding an appropriate plurisubharmonic function if necessary, we may assume that
$\phi$ reaches its minimum in an interior point $x\in B_0$.
Adding an appropriate constant, we can assume that
$U_0:= \phi^{-1}(-\infty, 0)$ has compact closure
with smooth strictly pseudoconvex boundary.
Let $U$ be the closure of $\Psi^{-1}(U_0)$,
and $\hat U$ the preimage of $U$ under the
bimeromorphic contraction $\hat B \arrow B$.
We obtain $M$ by gluing two boundary
components of $\hat U\backslash U_0$
as in \ref{_domain_C^2_Remark_}.

\hfill

{\bf Step 3:}
The map $\pi:\; \hat U \arrow U$ is a proper
holomorphic map of manifolds of the same dimension.
Therefore, the pushforward of a positive
$(p,p)$-form is a positive $(p,p)$-current.

Let $\pi_*\hat \omega$ be the pushforward of
$\hat \omega$, considered as a current on $U$. 
Using the $dd^c$-lemma for
currents, we obtain $\pi_*\hat\omega= dd^c f$, where
$f$ is a plurisubharmonic function on $U$ which is
smooth outside of the singularities of $\pi$.

Denote by $R$ the boundary of $U\subset \C^2$.
Then $f$ is smooth and strictly plurisubharmonic
in a neighbourhood of $R$.
Another plurisubharmonic function in a neighbourhood of
$R$ is obtained by taking $f_1:=(\Psi^{-1})^*\phi$,
where $\phi$ is the K\"ahler potential on
$U_0$ constructed in Step 2.

Rescaling $f$ if necessary and adding
a constant, we may assume that $-\epsilon <f\restrict R<0$
and $|df|\restrict R \ll \epsilon$. Let $A$ be a sufficiently
big positive number, and $0 < \delta \ll \epsilon$.
Then the regularized maximum $\max_{\delta}$ of
$f$ and $A f_1$
is equal to $A f_1$ in a very small neighbourhood of $R$
(because $f$ is negative on $R$ and $f_1=0$ on $R$),
and equal to $f$ in a neighbourhood $V$ of
$R_\epsilon:=Af_1^{-1}(-2\epsilon)$ because $|df| \ll A |df_1|$
and as $Af_1$ goes to $-2\epsilon$, $f$ does not go
below $-\epsilon$.

Replacing $\hat \omega$ by
$dd^c \max_{\delta}(f, f_1)$ on the annulus between $R$
and $R_{2\epsilon}$, we obtain a K\"ahler form $\hat \omega_1$.
Since $\max_{\delta}(f, f_1) = f_1$ in a neighbourhood of $R$,
the map $\Psi:\; (U, \hat\omega_1) \arrow (U_0, \hat \omega_1)$ acts
by homothety mapping a neighbourhood of $R$ with the
metric $\hat \omega_1$ isometrically to a neighbourhood of $\Psi(R)$
with the metric $A\hat \omega_1$ (however,
$A\hat \omega_1= A\hat \omega$ outside the annulus
bounded by $R$ and $R_{2\epsilon}$).
We have constructed an LCK metric on any Kato surface.
\endproof

\hfill

\remark 
The argument used in the proof of Brunella's theorem is
valid for any open ball in $\C^n$, giving a generalization
of Kato manifolds to arbitrary dimension. These manifolds,
were explored in \cite{_IOP:new_Kato_}. By Brunella's theorem,
they also admit an LCK structure.

\hfill

\remark Prior to publishing his general result
in \cite{_Brunella:Kato_}, M. Brunella constructed examples
of LCK metrics on Enoki surfaces
(\cite{_Brunella:Enoki_}), which
are special cases of class VII surfaces admitting 
a global spherical shell. Even before,  LCK metrics on some
Kato surfaces (hyperbolic Inoue and parabolic Inoue) also
appeared as anti-self-dual bihermitian metrics, obtained
by twistor methods in the work of A. Fujiki and
M. Pontecorvo, see  \cite{_Fujiki_Max:JDG_,
	_Fujiki_Max:SIGMA_}. In \cite{_Fujiki_Max:arxiv016_}
possible values of the Lee class (cohomology class of the
Lee form $\theta\in \Lambda^1(M)$) of these
LCK structures is discussed; see also
\cite{_Apostolov_Dloussky_}. As M. Pontecorvo pointed out, 
\cite{_Pontecorvo:letter_},
it is implicitly shown there that these Lee classes are different from the
ones constructed by Brunella. 

\hfill

\noindent{\bf Acknowledgments:}
We are grateful to Georges Dloussky and Matei Toma for a careful
reading of an early version of this paper and to Eduardo Esteves for an interesting discussion
about elliptic surfaces. Stefan Nemirovski and Iku Nakamura 
helped with reference and corrections to the first version
of this paper. Multiple thanks to Dmitry Kaledin and 
the anonymous referee for many insightful
comments and useful suggestions which 
improved many proofs.

\hfill

\hfill

{\small

\noindent {\sc Liviu Ornea\\
University of Bucharest, Faculty of Mathematics, \\14
Academiei str., 70109 Bucharest, Romania}, and:\\
{\sc Institute of Mathematics "Simion Stoilow" of the Romanian
Academy,\\
21, Calea Grivitei Str.
010702-Bucharest, Romania\\
\tt lornea@fmi.unibuc.ro, \ \  liviu.ornea@imar.ro}

\hfill

\noindent {\sc Misha Verbitsky\\
{\sc Instituto Nacional de Matem\'atica Pura e
              Aplicada (IMPA) \\ Estrada Dona Castorina, 110\\
Jardim Bot\^anico, CEP 22460-320\\
Rio de Janeiro, RJ - Brasil }\\
also:\\
Laboratory of Algebraic Geometry, \\
Faculty of Mathematics, National Research University 
Higher School of Economics,
7 Vavilova Str. Moscow, Russia}\\
\tt verbit@verbit.ru, verbit@impa.br }

\hfill

\noindent {\sc Victor Vuletescu\\
University of Bucharest, Faculty of Mathematics, \\14
Academiei str., 70109 Bucharest, Romania}\\
\tt vuli@fmi.unibuc.ro

\end{document}